\documentclass{amsart}
\usepackage{amsmath,amssymb, amsbsy}
\usepackage[dvips]{graphicx}
\newcommand{\R}{{\mathbb R}}
\newcommand{\N}{{\mathbb N}}

\renewcommand{\a }{\alpha }

\newcommand{\D }{\Delta }
\newcommand{\e }{\varepsilon }
\newcommand{\g }{\gamma}

\newcommand{\n }{\nabla }

\newcommand{\Di}{{\mathcal D}^{1,2}(\R^N)}
\newcommand{\alchi}{\raisebox{1.7pt}{$\chi$}}

\newenvironment{pf}{\noindent{\sc Proof}.\enspace}{\rule{2mm}{2mm}\medskip}
\newenvironment{pfn}[1]{\noindent{\bf Proof of {#1}.\enspace}}{\rule{2mm}{2mm}
\hfill\medskip}
\newtheorem{Theorem}{Theorem}[section]
\newtheorem{Corollary}[Theorem]{Corollary}

\newtheorem{Lemma}[Theorem]{Lemma}
\newtheorem{Proposition}[Theorem]{Proposition}

\newtheorem{remark}[Theorem]{Remark}
\begin{document}

\title[Nonlinear Schr\"odinger 
equations with inverse-square 
anisotropic 
potentials]{On the existence of ground state solutions to nonlinear Schr\"odinger 
equations with multisingular inverse-square 
anisotropic potentials}

\author[Veronica Felli]{Veronica Felli}

\address{Universit\`a degli Studi di Milano Bicocca, Dipartimento di Matematica e
  Applicazioni, Via Cozzi 53, 20125 Milano, Italy.}  \email{{\tt
    veronica.felli@unimib.it}. }

\date{February 1, 2008} 

\thanks{Supported by Italy MIUR, national project ``Variational
  Methods and Nonlinear Differential Equations''.
  }

\begin{abstract}
  \noindent A class of nonlinear
Schr\"odinger equations with  critical power-nonlinearities and 
potentials exhibiting multiple anisotropic inverse square
singularities is investigated. Conditions on 
strength, location, and orientation of singularities are given for  
the minimum of the associated  Rayleigh quotient to be 
achieved, both in the whole $\R^N$ and in bounded domains.
\end{abstract}
 \maketitle

\section{Introduction and statement of the main results}
This paper is concerned with the following class of nonlinear
Schr\"odinger equations with a critical power-nonlinearity and a
potential exhibiting multiple anisotropic inverse square
singularities:
\begin{equation}\label{eq:14}
\begin{cases}
  -\Delta
  v-{\displaystyle{\sum_{i=1}^k}}\dfrac{h_i\big(\frac{x-a_i}
{|x-a_i|}\big)}{|x-a_i|^2}\,v=v^{2^*-1},\\[15pt]
  v>0\quad\text{in }\R^N\setminus\{a_1,\dots,a_k\},
\end{cases}
\end{equation}
where $N\ge 3$,
 $k\in\N$, $h_i\in C^1({\mathbb S}^{N-1})$,
$(a_1,a_2,\dots, a_k)\in \R^{k N}$, $a_i\neq a_j$ for $i\neq j$, and
$2^*=\frac{2N}{N-2}$ is the critical Sobolev exponent.

The interest in such a class of equations arises in nonrelativistic
molecular physics.  Inverse square potentials with anisotropic
coupling terms turn out to describe the interaction between electric
charges and dipole moments of molecules, see \cite{leblond}. In
crystalline matter, the presence of many dipoles leads to consider
multisingular Schr\"odinger operators of the form
\begin{equation}\label{eq:op}
-\D
-{\displaystyle{\sum_{i=1}^k}}
\dfrac{\lambda_i\,(x-a_i)\cdot{\mathbf d}_i}{|x-a_i|^3},
\end{equation}
where $\lambda_i>0$, $i=1,\dots,k$, is proportional to the magnitude
of the $i$-th dipole and ${\mathbf d}_i$, $i=1,\dots,k$, is the unit
vector giving the orientation of the $i$-th dipole.

Schr\"odinger equations and operators with isotropic inverse-square
singular potentials have been largely investigated in the literature,
both in the case of one pole, see e.g. \cite{AFP,   GP, Jan,
  SM, terracini}, and in that of multiple singularities, see
\cite{caohan, chen, duyckaerts, esteban, FMT, FT}.  The
anisotropic case was first considered in \cite{terracini} where the
problem of existence of ground state solutions to (\ref{eq:14})  was discussed for $k=1$.
In \cite{FMT2}, an asymptotic formula for solutions to equation
associated with dipole-type Schr\"odinger operators near the
singularity was established. We also mention that positivity,
localization of binding and essential self-adjointness properties of a
class of Schr\"odinger operators with many anisotropic inverse square
singularities were investigated in \cite{FMT3}.

Ground state solutions to (\ref{eq:14}), i.e. solutions 
with the smallest energy, can be obtained through minimization of  the
associated Rayleigh quotient
\begin{equation}\label{eq:min}
S(h_1,h_2,\dots, h_k)=\inf_{u\in{\mathcal
    D}^{1,2}(\R^N)\setminus\{0\}}
\frac{{\mathcal Q}(u)
}{\big(\int_{\R^N}  |u|^{2^*}dx\big)^{2/2^*}},
\end{equation}
where $\Di$ denotes the closure space of
$C_0^{\infty}(\R^N)$ with respect to the norm
$$
\|u\|_{\Di}:=\bigg(\int_{\R^N}|\n u|^2\,dx\bigg)^{1/2},
$$
and ${\mathcal Q}:\Di\to\R$ is the quadratic form associated to 
the left-hand side of  equation \eqref{eq:14}, i.e.
\begin{equation}\label{eq:qf}
{\mathcal Q}(u):=\int_{\R^N} |\n  u|^2dx-\sum_{i=1}^k
  \int_{\R^N}\frac{h_i\big(\frac{x-a_i}
{|x-a_i|}\big)}{|x-a_i|^2}\,u^2(x)\,dx.
\end{equation}
Positive minimizers of (\ref{eq:min}) suitably rescaled give rise to
weak $\Di$-solutions to (\ref{eq:14}), which, by the Brezis-Kato Theorem
\cite{BrezisKato} 
and standard  elliptic regularity theory turn out be classical solutions in 
$\R^N\setminus\{a_1,\dots,a_k\}$.

The present paper means to extend to problems (\ref{eq:14}) and
(\ref{eq:min}) the analysis performed in \cite{FT} in the case of
locally isotropic inverse square potentials (i.e. for all $h_i$'s constant),
proving conditions on the strength, location and orientation of singularities 
for their solvability.

 A necessary condition for the
existence of positive classical solutions to \eqref{eq:14} in
$\R^N\setminus\{a_1,\dots,a_k\}$ is that ${\mathcal Q}$ is positive
semidefinite in $\Di$.
\begin{Proposition}\label{p:posde}
A necessary condition for the solvability of problem  (\ref{eq:14}) is
that the quadratic form
${\mathcal Q}(u)$ defined in \eqref{eq:qf}  
is positive semidefinite, i.e 
$$
{\mathcal Q}(u)\geq 0\quad\text{for all }u\in\Di.
$$
\end{Proposition}
A necessary  condition on the angular coefficients
$h_i$'s for the positive semidefiniteness of the quadratic form can be 
expressed in terms of the first eigenvalues of the associated Schr\"odinger 
operators on the sphere. Indeed, letting, 
for any $h\in C^1\big({\mathbb S}^{N-1}\big)$,   
$\mu_1(h)$ be the first eigenvalue of the
operator $-\D_{\mathbb S^{N-1}}-h(\theta)$ on~$\mathbb S^{N-1}$,
i.e. 
$$
\mu_1(h)=\min_{\psi\in H^1(\mathbb
S^{N-1})\setminus\{0\}}\frac{\int_{\mathbb S^{N-1}}|\n_{\mathbb
S^{N-1}}\psi(\theta)|^2\,dV(\theta)-\int_{\mathbb S^{N-1}}h(\theta)
\psi^2(\theta)\,dV(\theta)}{\int_{\mathbb S^{N-1}}\psi^2(\theta)\,dV(\theta)},
$$
a necessary (but not sufficient) condition for the quadratic form defined in (\ref{eq:qf}) 
to be  positive semidefinite is that
\begin{equation}\label{eq:30}
\mu_1(h_i)\geq-\left(\frac{N-2}2\right)^{\!\!2},\quad\text{for all
}i=1,\dots,k, \quad\text{and}\quad
\mu_1\left({\textstyle{\sum_{i=1}^k}
    h_i}\right)\geq-\left(\frac{N-2}2\right)^{\!\!2},
\end{equation}
see \cite{FMT2}. 
 
In particular, condition (\ref{eq:30}) is necessary for solvability 
of problem (\ref{eq:14}). In this paper, we shall actually consider 
multisingular anisotropic potentials with angular terms satisfying the
stronger assumption  
\begin{equation}\label{eq:31}
\mu_1(h_i)>-\left(\frac{N-2}2\right)^{\!\!2},\quad\text{for all
}i=1,\dots,k, \quad\text{and}\quad
\mu_1\left({\textstyle{\sum_{i=1}^k}
    h_i}\right)>-\left(\frac{N-2}2\right)^{\!\!2}.
\end{equation}
In \cite[Proposition 1.2]{FMT3} it was proved that condition
(\ref{eq:31}) is necessary for the quadratic form ${\mathcal Q}$ to be
positive definite, i.e. to have
\begin{equation}\label{eq:posde}
\mu(h_1,\dots,h_k,a_1,\dots,a_k):=
\inf_{\Di\setminus\{0\}}\frac{{\mathcal Q}(u)}{\|u\|_{\Di}^2}>0.
\end{equation}
On the other hand, (\ref{eq:31}) is not sufficient for the validity of
(\ref{eq:posde}), see \cite[Example 1.5]{FMT3}. However, if (\ref{eq:31})
holds, then (\ref{eq:posde}) turns out to be necessary for the solvability 
of (\ref{eq:14}).
\begin{Proposition}\label{p:posde2}
If (\ref{eq:31}) holds and (\ref{eq:14}) admits a positive $\Di$-solution, then 
(\ref{eq:posde}) is necessarily satisfied.
\end{Proposition}

\noindent Due to the above proposition, in order to look for solutions 
to (\ref{eq:14}), we will assume that 
the quadratic form ${\mathcal Q}$ is positive definite.
The dependence of positivity of the quadratic form on the
location and orientation of dipoles has been deeply investigated in
\cite{FMT3}, where conditions on the $h_i$'s and $a_i$'s ensuring the
validity of (\ref{eq:posde}) can be found. 
If ${\mathcal Q}(u)$ is
positive definite, then Sobolev's inequality implies that
$$
S(h_1,h_2,\dots, h_k)\geq
\mu(h_1,\dots,h_k,a_1,\dots,a_k) \,S>0,
$$
where $S$ is the best constant in the classical Sobolev inequality, i.e.
$$
S= \inf_{\Di\setminus\{0\}}\frac{\|u\|_{\Di}^2}{\|u\|_{L^{2^*}(\R^N)}^2}.
$$
Problems (\ref{eq:14}) and (\ref{eq:min}) have been treated by Terracini in 
\cite{terracini} in the one-dipole case $k=1$. For $h\in  C^1({\mathbb S}^{N-1})$, let
\begin{equation}\label{eq:2}
  S(h):=\inf_{u\in\Di\setminus\{0\}}\frac{\int_{\R^N}\big[
    |\nabla u(x)|^2-\frac{h(x/|x|)}{|x|^2}u^2(x)\big]\,dx}
  {\big(\int_{\R^N}|u|^{2^*}\big)^{2/2^*}}.
\end{equation}
Let us recall from \cite{terracini} the following existence result for the 
one-dipole type problem. 
\begin{Theorem}\label{l:t96}\cite[Proposition 5.3 and Theorem 0.2]{terracini}
Let $h\in C^1({\mathbb S}^N)$ such that $\mu_1(h)>-\big(\frac{N-2}2\big)^2$ and
\begin{equation}\label{eq:1}
\begin{cases}
\max_{{\mathbb S}^{N-1}}h>0,\quad&\text{if }N\geq 4,\\
\int_{{\mathbb S}^{N-1}}h\geq 0,\quad&\text{if }N=3.
\end{cases}
\end{equation}
Let $S(h)$ be defined in (\ref{eq:2}). Then $S(h)<S$ and $S(h)$ 
is achieved.
\end{Theorem}
The main difficulty in the minimization of the Rayleigh quotient in
(\ref{eq:min}) is due to the lack of compactness of the embeddings
$\Di\hookrightarrow L^{2^*}(\R^N)$ and $\Di\hookrightarrow
L^{2}\big(\R^N, |x|^{-2}h\big({x}/{|x|}\big)dx\big)$, where, for $h\in
L^{\infty}({\mathbb S}^{N-1})$, $L^{2}\big(\R^N,
|x|^{-2}h\big({x}/{|x|}\big)dx\big)$ is the the weighted Lebesgue
space endowed with the norm
$\big(\int_{\R^N}|x|^{-2}h\big({x}/{|x|}\big)u^2(x)\,dx\big)^{1/2}$.
Such a lack of compactness could produce non convergence of minimizing
sequences and non attainability of the infimum of the Rayleigh
quotient in some cases. In \cite{FT}
several configurations for which the infimum of the Rayleigh quotient is not attained 
are produced in the isotropic case, i.e. for all $h_i$'s constant;
e.g. the infimum in  (\ref{eq:min}) is not attained if the coefficients 
$h_i$'s are positive constants or if $k=2$ and $h_1$ and $h_2$ are costant.

A careful analysis of the behavior of minimizing sequences performed
through the P. L. Lions Concentration-Compactness Principle \cite{PL1,
  PL2} clarifies what are the possible reasons for lack of
compactness: concentration of mass at some non-singular point, at one
of the singularities or at infinity, see Theorem
\ref{th:ps}. Extending analogous results of \cite{FT} for the
isotropic case, Theorem \ref{t:ach} below provides sufficient
conditions for minimizing sequences to stay at an energy level which
is strictly below all the energy thresholds at which the compactness
can be lost. The proof is based on a comparison between levels which
is carried out by testing the energy functional associated to
(\ref{eq:14}) with solutions to (\ref{eq:2}).  On the other hand,
while in the isotropic case the solutions to (\ref{eq:2}) are
completely classified and can be explicitly written, in the
anisotropic case an explicit form of them is not available.  We
overcome this difficulty by exploiting the asymptotic analysis of the
behavior near the singularities of solutions performed in \cite{FMT2},
which allows us to estimate the behavior of minimizing sequences and
to force their level to stay in the recovered compactness range. 

From
now on, for every $h\in C^1({\mathbb S}^{N-1})$, we denote as
$\mu_1(h)$ the first eigenvalue of the operator $-\D_{\mathbb
  S^{N-1}}-h(\theta)$ on $\mathbb S^{N-1}$ and by $\psi_1^h$ the
associated positive $L^2$-normalized eigenfunction, and set
\begin{equation}\label{eq:sigmaH}
\sigma_h:=-\frac{N-2}2+\sqrt{\bigg(\frac{N-2}2\bigg)^{\!\!2}+\mu_1(h)}.
\end{equation}
\begin{Theorem}\label{t:ach}
  For $i=1,\dots,k$, let $a_i\in \R^{N}$, $a_i\neq a_j$ for $i\neq j$,  and 
  $h_i\in C^1({\mathbb S}^N)$  satisfy (\ref{eq:posde}).  If 
 \begin{align}
 \label{eq:23}
 &S(h_k)=\min \{S(h_j): j=1,\dots,k\},\\[7pt]
 \label{eq:24}
 &h_k \text{ satisfies (\ref{eq:1})},\\[7pt]
 &\label{eq:25}
 \begin{cases}
 {\displaystyle{\sum_{i=1}^{k-1}}}\,\dfrac{h_i\big(\frac{a_k-a_i}
 {|a_k-a_i|}\big)}{|a_k-a_i|^2}>0,&\text{if }
 \mu_1(h_k)\geq -\big(\frac{N-2}2\big)^2+1,\\[7pt]
 {\displaystyle{\sum_{i=1}^{k-1}\,
 \int_{\R^N}}} \dfrac{h_i\big(\frac{x}{|x|}\big)\Big[\psi_1^{h_k}
  \big(\frac{x+a_i-a_k}{|x+a_i-a_k|}\big)\Big]^2}
    {|x|^2|x+a_i-a_k|^{2(\sigma_{h_k}+N-2)}}>0,&\text{if }
 -\big(\frac{N-2}2\big)^2<\mu_1(h_k)< -\big(\frac{N-2}2\big)^2+1,
 \end{cases}\\[7pt]
&\label{eq:26}
S(h_k)\leq S\Big({\textstyle{\sum_{i=1}^k}} h_i  \Big),
 \end{align}
 then the infimum in (\ref{eq:min}) is achieved and problem
 (\ref{eq:14}) admits a solution in $\Di$.
\end{Theorem}
We notice that $S(h)=S(h\circ A)$ for any $h\in C^1({\mathbb S}^N)$
and any orthogonal matrix $A\in O(N)$. Hence condition \eqref{eq:26}
is satisfied for example if there exists an orthogonal matrix $A\in
O(N)$ such that
$$
\sum_{i=1}^k h_i(\theta)\leq h_k(A(\theta)),\quad\text{for all }\theta
\in{\mathbb S}^{N-1}.
$$
Let us describe in more detail  the case in which the singularities are generated
by electric dipoles, i.e. $h_i(\theta)=\lambda_i\theta\cdot
{\mathbf d}_i$, for some $\lambda_i>0$ and ${\mathbf d}_i\in\R^N$ with
$|{\mathbf d}_i|=1$. For any $\lambda>0$ and ${\mathbf d}\in\R^N$ with
$|{\mathbf d}|=1$, let  
\begin{equation*}\label{firsteig}
  \mu_1^\lambda=\min_{\psi\in H^1(\mathbb
    S^{N-1})\setminus\{0\}}\frac{\int_{\mathbb S^{N-1}}|\n_{\mathbb
      S^{N-1}}\psi(\theta)|^2\,dV(\theta)-\lambda \int_{\mathbb S^{N-1}}(\theta\cdot
    {\mathbf d})\psi^2(\theta)\,dV(\theta)}
  {\int_{\mathbb S^{N-1}}\psi^2(\theta)\,dV(\theta)}
\end{equation*}
be the first eigenvalue of the
operator $-\D_{\mathbb S^{N-1}}-\lambda\,(\theta\cdot {\mathbf d})$ on $\mathbb S^{N-1}$.
By rotation invariance, it is easy to verify that the above minimum  
does not depend on ${\mathbf d}$. Moreover, condition (\ref{eq:31}) can be 
explicitly expressed as a bound on the dipole magnitudes; indeed, 
$$
\mu_1^{\lambda}>-\left(\frac{N-2}2\right)^{\!\!2}
\quad\text{if and only if}\quad
\lambda<\frac1{\Lambda_N}
$$
where $\Lambda_N$ is  the best constant 
in the dipole Hardy-type inequality, i.e.
$$
\Lambda_N:=\sup_{u\in\Di\setminus\{0\}}\dfrac{{\displaystyle
{\int_{\R^N}{\dfrac{x\cdot{\mathbf
d}}{|x|^3}\,u^2(x)\,dx}}}}{{\displaystyle{\int_{\R^N}{|\n u(x)|^2\,dx}}}},
$$
see \cite{FMT2}.  By rotation invariance,
$\Lambda_N$ does not depend on the unit vector ${\mathbf d}$ and, by
classical Hardy's inequality, $\Lambda_N<4/(N-2)^2$. For every
$\lambda>0$, let us denote
$\sigma^\lambda:=-\frac{N-2}2+\sqrt{\big(\frac{N-2}2\big)^{2}+\mu_1^{\lambda}}$.

\begin{Corollary}\label{c:ach}
  For $i=1,\dots,k$, let $a_i\in \R^{N}$, $a_i\neq a_j$ for $i\neq j$,
${\mathbf d}_i\in \R^{N}$ with 
$|{\mathbf d}_i|=1$, 
 and
$$
0<\lambda_1\leq \lambda_2\leq\dots\leq\lambda_k<\Lambda_N^{-1}.
$$
Assume that the quadratic form 
$$
u\mapsto\int_{\R^N} |\n
u(x)|^2dx-{\displaystyle{\sum_{i=1}^k}}
\dfrac{\lambda_i\,(x-a_i)\cdot{\mathbf d}_i}{|x-a_i|^3}\,u^2(x)\,dx
$$ 
is positive definite and that
 \begin{align}
 &\label{eq:27}
 \begin{cases}
 {\displaystyle{\sum_{i=1}^{k-1}}}\,\dfrac{\lambda_i\,{\mathbf d}_i\cdot
\frac{a_k-a_i}
 {|a_k-a_i|}}{|a_k-a_i|^2}>0,&\text{if }
 \mu_1^{\lambda_k}\geq -\big(\frac{N-2}2\big)^2+1,\\[7pt]
 {\displaystyle{\sum_{i=1}^{k-1}\,
 \int_{\R^N}}} \dfrac{\lambda_i\frac{x}{|x|}\cdot{\mathbf d}_i
\Big[\psi_1^{\lambda_k \theta\cdot {\mathbf d}_k}
  \big(\frac{x+a_i-a_k}{|x+a_i-a_k|}\big)\Big]^2}
    {|x|^2|x+a_i-a_k|^{2(\sigma^{\lambda_k}+N-2)}}>0,&\text{if }
 -\big(\frac{N-2}2\big)^2<\mu_1^{\lambda_k}< -\big(\frac{N-2}2\big)^2+1,
 \end{cases}\\[7pt]
&\label{eq:29}
\left|\sum_{i=1}^k\lambda_i{\mathbf d}_i\right|\leq \lambda_k.
 \end{align}
Then
 the infimum 
$$
\inf_{u\in{\mathcal
    D}^{1,2}(\R^N)\setminus\{0\}}
\frac{\int_{\R^N} |\n
u(x)|^2dx-{\displaystyle{\sum_{i=1}^k}}
\dfrac{\lambda_i\,(x-a_i)\cdot{\mathbf d}_i}{|x-a_i|^3}\,u^2(x)\,dx
}{\big(\int_{\R^N}  |u|^{2^*}dx\big)^{2/2^*}}
$$
 is achieved and  the problem
\begin{equation}\label{eq:28}
\begin{cases}
-\D u
-{\displaystyle{\sum_{i=1}^k}}
\dfrac{\lambda_i\,(x-a_i)\cdot{\mathbf d}_i}{|x-a_i|^3}\,u
=u^{2^*-1},\\[15pt]
  u>0\quad\text{in }\R^N\setminus\{a_1,\dots,a_k\},
\end{cases}
\end{equation}
admits a solution in $\Di$.
\end{Corollary}
With respect to the isotropic case, the possibility of orientating the
dipoles helps in finding the balance between the strength and the
locations of the singularities required in assumptions
(\ref{eq:27}--\ref{eq:29}).  Let us consider for example the case of
two dipoles $k=2$.  Assume that $0<\lambda_1\leq\lambda_2$,
$\lambda_2$ is small and $N$ is large in such a way that the
associated quadratic form is positive definite and
$\mu_1^{\lambda_2}\geq-\big(\frac{N-2}2\big)^2+1$. Then condition
(\ref{eq:27}) reads as
 $$
 (a_2-a_1)\cdot {\mathbf d}_1>0,
  $$
while (\ref{eq:29}) reads as
 $$
 {\mathbf d}_1\cdot{\mathbf d}_2<-\frac{\lambda_1}{2\lambda_2}.
 $$
 In this case, if the first dipole $\lambda_1 {\mathbf d}_1$ is fixed
 at point $a_1$, (\ref{eq:27}) gives a constraint on the location of
 the second dipole while (\ref{eq:29}) gives a condition on its
 orientation. In particular, it is possible to construct many
 configurations ensuring the existence of ground state solutions to
 (\ref{eq:28}), unlike the isotropic case where problem (\ref{eq:14})
 with $k=2$ and $h_1$ and $h_2$ constants has no ground state
 solutions, as observed in \cite[Theorem 1.3]{FT}.

In  bounded domains,   concentration of mass at  infinity
is no more possible and an existence result similar to Theorem \ref{t:ach}
can be obtained without assumption (\ref{eq:26}). 
\begin{Theorem}\label{t:exlim}
Assume that  $\Omega$ is a bounded smooth domain, 
$\{a_1,a_2,\dots,
a_k\}\subset\Omega$,  
   $h_i\in C^1({\mathbb S}^N)$,  $i=1,\dots,k$, 
such that 
the quadratic form 
\begin{equation}\label{eq:pdbd}
{\mathcal Q}_{\Omega}(u)=
:=\int_{\Omega} |\n  u(x)|^2dx-
\sum_{i=1}^k
\int_{\R^N}
\dfrac{h_i\big(\frac{x-a_i}
{|x-a_i|}\big)}{|x-a_i|^2}u^2(x)\,dx
\quad\text{is positive definite},
\end{equation}
$h_k$ satisfies (\ref{eq:1}), $S(h_k)=\min \{S(h_j): j=1,\dots,k\}$,
$$
 \mu_1(h_k)\geq -\bigg(\frac{N-2}2\bigg)^{\!\!2}+1,\quad\text{and}\quad
\sum_{i=1}^{k-1}\,\dfrac{h_i\big(\frac{a_k-a_i}
 {|a_k-a_i|}\big)}{|a_k-a_i|^2}>0.
$$
Then the infimum in
\begin{equation}\label{eq:minbound}
S_{\Omega}(h_1,h_2,\dots, h_k)=\inf_{u\in H^1_0(\Omega)\setminus\{0\}}
\frac{{\mathcal Q}_{\Omega}(u)}{\|u\|_{L^{2^*}(\Omega)}^2},
\end{equation}
 is achieved and  equation
 \begin{equation}\label{eq:bound}
\begin{cases}
-\D
u-
{\displaystyle{\sum_{i=1}^k}}
\dfrac{h_i\big(\frac{x-a_i}
{|x-a_i|}\big)}{|x-a_i|^2}\,u=u^{2^*-1},\\[15pt]
u>0\quad\text{in }\Omega\setminus\{a_1,\dots,a_k\},\quad
u=0\quad\text{on } \partial\Omega,
\end{cases}
\end{equation}
 admits a solution in $H^1_0(\Omega)$.
\end{Theorem}

The further assumption  $\mu_1(h_k)\geq -\big(\frac{N-2}2\big)^{2}+1$ of 
Theorem \ref{t:exlim} is not technical but quite natural when working in 
bounded domains. Indeed it plays the role of a critical dimension for 
Brezis-Nirenberg type problems in bounded domains, see  \cite{BN, Jan}.

The paper is organized as follows. Section \ref{sec:necess-posit-quadr} 
contains the proofs of Propositions \ref{p:posde} and \ref{p:posde2}. In 
section \ref{sec:inter-estim} some interaction estimates are first deduced 
and then applied to comparison of energy levels of minimizing sequences. 
Section \ref{sec:pala-smale-cond} provides  a local Palais-Smale condition
which is  used to prove Theorem  \ref{t:ach} and Corollary \ref{c:ach}.
Finally, in section \ref{sec:probl-bound-doma} we analyze the problem in bounded 
domains.

\medskip
\noindent
{\bf Notation. } We list below some notation used throughout the
paper.\par
\begin{itemize}
\item[-]$B(a,r)$ denotes the ball $\{x\in\R^N: |x-a|<r\}$ in $\R^N$ with
center at $a$ and radius $r$.
\item[-] For any $A\subset \R^N$, $\alchi_{A}$ denotes the
characteristic function of $A$. 
\item[-] $S$ is the best constant in the Sobolev inequality
$S\|u\|_{L^{2^*}(\R^N)}^2\leq \|u\|_{{\mathcal D}^{1,2}(\R^N)}^2$.
\item[-] $\omega_N$ denotes the volume of the unit ball in $\R^N$.
\item[-] $O(N)$ denotes the group of orthogonal $N\times N$ matrices.
\end{itemize}

\section{Necessity of the positivity of the quadratic form}
\label{sec:necess-posit-quadr}
In the present section we discuss the necessity of the positivity of
the quadratic form for the solvability of (\ref{eq:14}), by proving
Propositions \ref{p:posde} and \ref{p:posde2}.

\begin{pfn}{Proposition \ref{p:posde}}
Let $u$ be a  positive classical solutions to \eqref{eq:14} in
$\R^N\setminus\{a_1,\dots,a_k\}$. For any $\phi\in C^{\infty}_c(\R^N\setminus
\{a_1,\dots,a_k\})$, by testing equation
\eqref{eq:1} with  $\frac{\phi^2}u$ we obtain
$$
2\int_{\R^N}\frac{\phi}u\,\nabla\phi\cdot\nabla
u\,dx-\int_{\R^N}\frac{\phi^2}{u^2}\,|\nabla
u|^2\,dx-\sum_{i=1}^k
  \int_{\R^N}\frac{h_i\big(\frac{x-a_i}
{|x-a_i|}\big)}{|x-a_i|^2}\,u^2(x)\,dx
-\int_{\R^N}\phi^2u^{2^*-2}\,dx=0.
$$
From the elementary inequality $2\,\frac{\phi}u\,\nabla\phi\cdot\nabla
u-\frac{\phi^2}{u^2}\,|\nabla u|^2\leq |\nabla\phi|^2$, we deduce
$$
{\mathcal Q}(\phi)\geq\int_{\R^N}\phi^2u^{2^*-2}\,dx\geq 0\quad\text{for all }\phi\in
C_c^{\infty}(\R^N\setminus 
\{a_1,\dots,a_k\}) .
$$
From density of $C_c^{\infty}(\R^N\setminus
\{a_1,\dots,a_k\})$  in $\Di$ (see \cite[Lemma
2.1]{catrinawang}), we obtain that  ${\mathcal Q}$ is positive semidefinite 
in $\Di$.
\end{pfn}

\begin{pfn}{Proposition \ref{p:posde2}}
  Assume that (\ref{eq:31}) holds and let $u\in\Di$ be a positive
  $\Di$-solution to (\ref{eq:14}). By Proposition \ref{p:posde}, it follows that
$$
\mu(h_1,\dots,h_k,a_1,\dots,a_k)\geq 0,
$$
where $\mu(h_1,\dots,h_k,a_1,\dots,a_k)$ has been defined in (\ref{eq:posde}). 
Let us assume by contradiction that $\mu(h_1,\dots,h_k,a_1,\dots,a_k)=0$. 
From \cite[Proposition 4.1]{FMT3}, $\mu(h_1,\dots,h_k,a_1,\dots,a_k)=0$ is attained by 
some $v\in\Di$, $v\geq 0$ a.e. in $\R^N$, $v\not\equiv 0$, which then satisfies
$$
-\Delta v-{\displaystyle{\sum_{i=1}^k}}\dfrac{h_i\big(\frac{x-a_i}
{|x-a_i|}\big)}{|x-a_i|^2}\,v=0\quad\text{weakly in }\Di.
$$
Testing the above equation with $u$, we obtain that 
$$
\int_{\R^N}
u^{2^*-1}(x)v(x)\,dx=0
$$
which is in contradiction with the positivity of $u$.
\end{pfn}
\section{Interaction estimates and comparison of energy
  levels}\label{sec:inter-estim}

By Theorem \ref{l:t96}, for every function $h\in C^1({\mathbb S}^N)$ verifying
$\mu_1(h)>-\big(\frac{N-2}2\big)^2$
 and (\ref{eq:1}), there exists some $\phi_h\in\Di$,
$\phi_h\geq0$, $\phi_h\not\equiv 0$, such that $\phi_h$ attains
$S(h)$, i.e.
\begin{equation}\label{eq:17bi}
S(h)=\frac{\int_{\R^N}\big[
|\nabla \phi_h(x)|^2-\frac{h(x/|x|)}{|x|^2}\phi_h^2(x)\big]\,dx}
{\big(\int_{\R^N}|\phi_h|^{2^*}\big)^{2/2^*}},
\end{equation}
and solves
\begin{equation}\label{eq:18bi}
-\Delta\phi_h-\frac{h(x/|x|)}{|x|^2}\phi_h=\phi_h^{2^*-1},\quad
\text{in }\R^N.
\end{equation}
Moreover the Kelvin's transform $w_h(x):=|x|^{-(N-2)}\phi_h(x/|x|^2)$
solves $-\Delta w_h-\frac{h(x/|x|)}{|x|^2}w_h=w_h^{2^*-1}$
in $\R^N$.
From \cite{FMT2}, it follows that, letting $\sigma_h$ defined in (\ref{eq:sigmaH}),
the functions 
$$
x\mapsto \frac{\phi_h(x)}{|x|^{\sigma_h}\psi_1^h(x/|x|)},\quad
x\mapsto \frac{w_h(x)}{|x|^{\sigma_h}\psi_1^h(x/|x|)}=
\frac{\phi_h(x/|x|^2)}{|x|^{\sigma_h+N-2}\psi_1^h(x/|x|)}
$$
are  continuous in $\R^N$ and admit positive limits as $|x|\to 0$,
i.e. 
\begin{equation}\label{eq:7bi}
c_0^h :=\lim_{|x|\to 0}
\frac{\phi_h(x)}{|x|^{\sigma_h}\psi_1^h(x/|x|)}\in(0,+\infty)
 \ \text{and}\ 
c_\infty^h:=\lim_{|x|\to +\infty}
\frac{\phi_h(x)}{|x|^{-\sigma_h-N+2}\psi_1^h(x/|x|)}\in(0,+\infty).
\end{equation}
Hence there exists a positive constant $C(h)>0$ such that
\begin{equation}\label{eq:3}
\frac1{C(h)}\,\frac{|x|^{\sigma_h}}{1+|x|^{2\sigma_h+N-2}}
\leq \phi_h(x)\leq \frac{C(h)\,|x|^{\sigma_h}}{1+|x|^{2\sigma_h+N-2}}, \quad
\text{for all }x\in\R^N\setminus\{0\}.
\end{equation}
For any $\mu>0$, let us denote $\phi_{\mu}^{h}(x):=\mu^{-(N-2)/2}\phi_h(x/\mu)$.
\begin{Lemma}\label{l:interest1}
  Let $h,k\in C^1({\mathbb S}^N)$ such that $h$ satisfies (\ref{eq:1})
  and $\mu_1(h)>-\big(\frac{N-2}2\big)^2+1$.  Then, for every
  $a\in\R^N\setminus\{0\}$, 
there holds
$$
\int_{\R^N}\phi_h^2(x)\,dx\in(0,+\infty)
\quad\text{and}\quad\int_{\R^N}\frac{k(\frac{x-a}{|x-a|})}{|x-a|^2}
|\phi_{\mu}^{h}(x)|^2\,dx=
\mu^2\left[
\frac{k\big(\frac{-a}{|a|}\big)}{|a|^2}\int_{\R^N}\phi_h^2(x)\,dx+o(1)\right]
$$
as $\mu\to 0^+$.
\end{Lemma}

\begin{pf}
From \eqref{eq:3} and the assumption $\mu_1(h)>-\big(\frac{N-2}2\big)^2+1$,
it follows that $\phi_h\in L^2(\R^N)$.
We have that 
\begin{align}\label{eq:22bi}
  \int_{\R^N}\frac{k(\frac{x-a}{|x-a|})}{|x-a|^2}
  |\phi_{\mu}^{h}(x)|^2\,dx
  &=\mu^2\int_{|x|<\frac{|a|}{2\mu}}\frac{k(\frac{\mu x-a}{|\mu
      x-a|})}{|\mu x-a|^2}\phi_h^2(x)\,dx \\
&\quad+\mu^{-N+2}\int_{|x+a|\geq\frac{|a|}{2}}
  \frac{k(\frac{x}{|x|})}{|x|^2}\phi_h^2\Big(\frac{x+a}{\mu}\Big)
\,dx
  \notag.
\end{align}
Since 
$$
\Bigg|\alchi_{B\big(0,\frac{|a|}{2\mu}\big)}(x)\frac{k\big(\frac{\mu x-a}{|\mu
      x-a|}\big)}{|\mu x-a|^{2}}\Bigg|\leq 
\frac4{|a|^{2}}\|k\|_{L^{\infty}({\mathbb S}^{N-1})}
$$
 and $\phi_h\in
 L^2(\R^N)$, from the Dominated Convergence Theorem we deduce that
\begin{align}\label{eq:4bi}
  \lim_{\mu\to 0^+}\int_{|x|<\frac{|a|}{2\mu}}\frac{k(\frac{\mu
      x-a}{|\mu x-a|})}{|\mu x-a|^2}\phi_h^2(x)\,dx &=\lim_{\mu\to
    0^+}\int_{\R^N}\alchi_{B\big(0,\frac{|a|}{2\mu}\big)}(x)\frac{k(\frac{\mu
      x-a}{|\mu
      x-a|})}{|\mu x-a|^2}\phi_h^2(x)\,dx\\
  &\notag=
  \frac{k\big(\frac{-a}{|a|}\big)}{|a|^2}\int_{\R^N}\phi_h^2(x)\,dx.
\end{align}
Moreover, from \eqref{eq:3} and  $\mu_1(h)>-\big(\frac{N-2}2\big)^2+1$, 
it follows that
\begin{multline}\label{eq:5bi}
\bigg|\mu^{-N}\int_{|x+a|\geq\frac{|a|}{2}}
  \frac{k(\frac{x}{|x|})}{|x|^2}\phi_h^2\Big(\frac{x+a}{\mu}\Big)
\,dx\bigg|\\
\leq \mu^{2\sigma_h+N-4}\|k\|_{L^{\infty}({\mathbb S}^{N-1})}(C(h))^2
\int_{|x+a|\geq\frac{|a|}{2}}
  \frac{1}{|x|^2 |x+a|^{2(\sigma_h+N-2)}}
\,dx=o(1)\quad\text{as }\mu\to 0^+.
\end{multline}
The conclusion follows then from 
\eqref{eq:22bi}, \eqref{eq:4bi}, and \eqref{eq:5bi}.
\end{pf}

\begin{Lemma}\label{l:interest2}
Let $h,k\in C^1({\mathbb S}^N)$ such that $h$ satisfies 
(\ref{eq:1}), and $\mu_1(h)=-\big(\frac{N-2}2\big)^2+1$. 
Then, for every $a\in\R^N\setminus\{0\}$, 
\begin{equation}\label{eq:6bi}
  N\omega_N\frac{|\log\mu|}{|C(h)|^2}\,(1+o(1))
  \leq\int_{|x|<\frac{|a|}{2\mu}}\phi_h^2(x)\,dx
  \leq N\omega_N|C(h)|^2\,|\log\mu|\,(1+o(1)),\quad
\text{as }\mu\to 0^+,
\end{equation}
where $\omega_N$ is the volume of the standard unit $N$-ball, and
\begin{equation}\label{eq:8bi}
\int_{\R^N}\frac{k(\frac{x-a}{|x-a|})}{|x-a|^2}
|\phi_{\mu}^{h}(x)|^2\,dx=
\mu^2
\bigg(\frac{k\big(\frac{-a}{|a|}\big)}{|a|^2}+o(1)\bigg)\bigg[\int_{|x|<\frac{1}{\mu}}
\phi_h^2(x)\,dx\bigg]
\end{equation}
as $\mu\to 0^+$.
\end{Lemma}
\begin{pf}
Estimate (\ref{eq:6bi}) follows from (\ref{eq:3}) and direct calculations.
We have that 
\begin{multline}\label{eq:22}
   \int_{\R^N}\frac{k(\frac{x-a}{|x-a|})}{|x-a|^2}
   |\phi_{\mu}^{h}(x)|^2\,dx\\
   =\mu^2\bigg[
 \frac{k\big(\frac{-a}{|a|}\big)}{|a|^2}
 \int_{|x|<\frac{|a|}{2\mu}}\phi_h^2(x)\,dx 
 +
 \int_{|x|<\frac{|a|}{2\mu}}k\Big(\frac{\mu x-a}{|\mu
       x-a|}\Big)\bigg(\frac{1}{|\mu x-a|^2}-\frac{1}{|a|^2}\bigg)\phi_h^2(x)\,dx 
\\ +
 \frac1{|a|^2}
 \int_{|x|<\frac{|a|}{2\mu}}\bigg(k\Big(\frac{\mu x-a}{|\mu
       x-a|}\Big)-k\Big(\frac{-a}{|a|}\Big)\bigg)
 \phi_h^2(x)\,dx 
 +\mu^{-N}\int_{|x+a|\geq\frac{|a|}{2}}
   \frac{k(\frac{x}{|x|})}{|x|^2}\phi_h^2\Big(\frac{x+a}{\mu}\Big)
 \,dx\bigg].
\end{multline}
Since 
$$
\bigg|\frac{1}{|\mu x-a|^2}-\frac{1}{|a|^2}\bigg|\leq
\frac{4}{|a|^4}\big(\mu^2|x|^2+2\mu|a||x|\big)
\quad\text{for }|x|<\frac{|a|}{2\mu},
$$
from \eqref{eq:3} it follows that
\begin{equation}\label{eq:9bi}
 \int_{|x|<\frac{|a|}{2\mu}}k\Big(\frac{\mu x-a}{|\mu
       x-a|}\Big)\bigg(\frac{1}{|\mu x-a|^2}-\frac{1}{|a|^2}\bigg)\phi_h^2(x)\,dx 
=O(1)\quad\text{as }\mu\to 0^+.
\end{equation}
Since $k\in C^1({\mathbb S}^N)$, for some positive constant $C$ depending on $k$
there holds
\begin{align*}
\bigg|k\Big(\frac{\mu x-a}{|\mu
       x-a|}\Big)-k\Big(\frac{-a}{|a|}\Big)\bigg|&\leq
C\, \bigg|\frac{\mu x-a}{|\mu
       x-a|}-\frac{-a}{|a|}\bigg|
= \frac{C\,\sqrt2}{\sqrt{|\mu
       x-a|}}\sqrt{|\mu
       x-a|-|a|+\mu\,\frac{a\cdot x}{|a|}}\\
&\leq \frac{C\,\sqrt2}{\sqrt{|\mu
       x-a|}}\sqrt{2\mu|x|}\leq \frac{2\,C\sqrt2\,\sqrt\mu\sqrt{|x|}}{\sqrt{|a|}}
\quad\text{for }|x|<\frac{|a|}{2\mu},
   \end{align*}
hence, from \eqref{eq:3}, it follows that
\begin{equation}\label{eq:10bi}
\int_{|x|<\frac{|a|}{2\mu}}\bigg(k\Big(\frac{\mu x-a}{|\mu
       x-a|}\Big)-k\Big(\frac{-a}{|a|}\Big)\bigg)
 \phi_h^2(x)\,dx =O(1)\quad\text{as }\mu\to 0^+.
\end{equation}
From \eqref{eq:3}, we deduce that
\begin{align*}
  \bigg|\mu^{-N}\int_{|x+a|\geq\frac{|a|}{2}}
  \frac{k(\frac{x}{|x|})}{|x|^2}\phi_h^2\Big(\frac{x+a}{\mu}\Big)
  \,dx\bigg| \leq\|k\|_{L^{\infty}({\mathbb S}^{N-1})}(C(h))^2
  \int_{|x+a|\geq\frac{|a|}{2}} \frac{1}{|x|^2 |x+a|^{N}}
\end{align*}
hence
\begin{equation}\label{eq:11bi}
  \mu^{-N}\int_{|x+a|\geq\frac{|a|}{2}}
  \frac{k(\frac{x}{|x|})}{|x|^2}\phi_h^2\Big(\frac{x+a}{\mu}\Big)
  \,dx=O(1)\quad\text{as }\mu\to 0^+.
\end{equation}
From 
\eqref{eq:6bi}, \eqref{eq:22}, \eqref{eq:9bi}, \eqref{eq:10bi}, and \eqref{eq:11bi} it
follows that
\begin{equation}\label{eq:19bi}
  \int_{\R^N}\frac{k(\frac{x-a}{|x-a|})}{|x-a|^2}
  |\phi_{\mu}^{h}(x)|^2\,dx=
  \mu^2
  \bigg(\frac{k\big(\frac{-a}{|a|}\big)}{|a|^2}+o(1)\bigg)
  \bigg[\int_{|x|<\frac{|a|}{2\mu}}
  \phi_h^2(x)\,dx\bigg]
\end{equation}
as $\mu\to 0^+$. From (\ref{eq:3}) and the assumption 
 $\mu_1(h)=-\big(\frac{N-2}2\big)^2+1$, we obtain that
$$
\left|\int_{|x|<\frac{|a|}{2\mu}}
  \phi_h^2(x)\,dx-\int_{|x|<\frac{1}{\mu}} \phi_h^2(x)\,dx\right|\leq
N\omega_N(C(h))^2\left|\int_{\frac{1}{\mu}}^{\frac{|a|}{2\mu}}r^{-1}dr\right|
=N\omega_N(C(h))^2\bigg|\log\frac{|a|}2\bigg|,
$$
hence, taking into account that, under assumption
$\mu_1(h)=-\big(\frac{N-2}2\big)^2+1$, $\phi_h\not\in L^2(\R^N)$,
\begin{equation}\label{eq:20}
  \int_{|x|<\frac{|a|}{2\mu}}
  \phi_h^2(x)\,dx=\int_{|x|<\frac{1}{\mu}}
  \phi_h^2(x)\,dx+O(1)=(1+o(1))\int_{|x|<\frac{1}{\mu}}
  \phi_h^2(x)\,dx
\end{equation}
as $\mu\to 0^+$. The conclusion (\ref{eq:8bi}) follows from (\ref{eq:19bi}) and 
(\ref{eq:20}).
\end{pf}

\begin{Lemma}\label{l:interest3}
  Let $h,k\in C^1({\mathbb S}^N)$ such that $h$ satisfies (\ref{eq:1})
  and $-\big(\frac{N-2}2\big)^2<\mu_1(h)<-\big(\frac{N-2}2\big)^2+1$.
  Then, for every $a\in\R^N\setminus\{0\}$ and $A\in O(N)$ such that
  $A e_1=\frac{a}{|a|}$, with $e_1=(1,0,\dots,0)\in\R^{N}$, there
  holds
\begin{align*}
  \int_{\R^N}\frac{k(\frac{x-a}{|x-a|})}{|x-a|^2}
  |\phi_{\mu}^{h}(x)|^2\,dx&= \mu^{2\sigma_h+N-2}\left[
    (c_\infty^h)^2\int_{\R^N}
    \frac{k\big(\frac{x}{|x|}\big)\Big[\psi_1^h
      \big(\frac{x+a}{|x+a|}\big)\Big]^2}
    {|x|^2|x+a|^{2(\sigma_h+N-2)}}
    \,dx+o(1)\right]\\
  &= \frac{\mu^{2\sigma_h+N-2}}{|a|^{2\sigma_h+N-2}} \left[
    (c_\infty^h)^2\int_{\R^N} \frac{(k\circ
      A)\big(\frac{x}{|x|}\big)\Big[(\psi_1^h\circ
      A)\big(\frac{x+e_1}{|x+e_1|}\big)\Big]^2}
    {|x|^2|x+e_1|^{2(\sigma_h+N-2)}}\,dx+o(1)\right]
\end{align*}
as $\mu\to 0^+$, being $c_\infty^h$ defined by (\ref{eq:7bi}).
\end{Lemma}
\begin{pf}
A direct calculation yields
\begin{align}\label{eq:12}
  \int_{\R^N}\frac{k(\frac{x-a}{|x-a|})}{|x-a|^2}
  |\phi_{\mu}^{h}(x)|^2\,dx& =\mu^{2\sigma_h+N-2}\int_{\R^N}
  \frac{k\big(\frac{x}{|x|}\big)\Big[\psi_1^h
    \big(\frac{x+a}{|x+a|}\big)\Big]^2} {|x|^2|x+a|^{2(\sigma_h+N-2)}}
  \cdot \frac{\phi_h^2\big(\frac{x+a}\mu\big)}
  {\big|\frac{x+a}\mu\big|^{2(2-\sigma_h-N)}\Big[\psi_1^h
    \big(\frac{x+a}{|x+a|}\big)\Big]^2} \,dx.
\end{align}
From \eqref{eq:3}, it follows that the function
$$
x\mapsto \frac{\phi_h^2\big(\frac{x+a}\mu\big)}
{\big|\frac{x+a}\mu\big|^{2(2-\sigma_h-N)}\Big[\psi_1^h
  \big(\frac{x+a}{|x+a|}\big)\Big]^2}
$$
is bounded a.e. in $\R^N$ uniformly with respect to $\mu>0$,
whereas \eqref{eq:7bi}  implies that, for a.e. $x\in\R^N$,
\begin{equation}\label{eq:13}
  \lim_{\mu\to0}\frac{\phi_h^2\big(\frac{x+a}\mu\big)}
  {\big|\frac{x+a}\mu\big|^{2(2-\sigma_h-N)}\Big[\psi_1^h
    \big(\frac{x+a}{|x+a|}\big)\Big]^2}=
  (c_\infty^h)^2 .
\end{equation}
Since the assumption $\mu_1(h)<-\big(\frac{N-2}2\big)^2+1$
ensures that 
$$
x\mapsto\frac{k\big(\frac{x}{|x|}\big)\Big[\psi_1^h
  \big(\frac{x+a}{|x+a|}\big)\Big]^2}
{|x|^2|x+a|^{2(\sigma_h+N-2)}}\in L^1(\R^N),
$$
from \eqref{eq:12}, \eqref{eq:13}, and the Dominated Convergence Theorem we deduce that
\begin{align*}
  \int_{\R^N}\frac{k(\frac{x-a}{|x-a|})}{|x-a|^2}
  |\phi_{\mu}^{h}(x)|^2\,dx& =\mu^{2\sigma_h+N-2}\left[
    (c_\infty^h)^2\int_{\R^N}
    \frac{k\big(\frac{x}{|x|}\big)\Big[\psi_1^h
      \big(\frac{x+a}{|x+a|}\big)\Big]^2}
    {|x|^2|x+a|^{2(\sigma_h+N-2)}} \,dx+o(1)\right]
\end{align*}
as $\mu\to0$.
Through the change of variable $x=|a|Ay$, we obtain that
\begin{align*}
  \int_{\R^N} \frac{k\big(\frac{x}{|x|}\big)\Big[\psi_1^h
    \big(\frac{x+a}{|x+a|}\big)\Big]^2} {|x|^2|x+a|^{2(\sigma_h+N-2)}}
  \,dx=|a|^{-N-2\sigma_h+2} \int_{\R^N} \frac{(k\circ
    A)\big(\frac{y}{|y|}\big)\Big[\psi_1^h \big(A
    \big(\frac{y+e_1}{|y+e_1|}\big)\big)\Big]^2}
  {|y|^2|y+e_1|^{2(\sigma_h+N-2)}} \,dx
\end{align*}
thus completing the proof.
\end{pf}

The interaction estimates provided by Lemmas \ref{l:interest1},
\ref{l:interest2}, and  \ref{l:interest3}, allow us to compare 
the ground state level of the multisingular problem   with 
the ground state level of the single dipole  problem.
\begin{Proposition}\label{cor:sub}
Let $h_i\in C^1({\mathbb S}^N)$, $i=1,\dots,k$.
Let us assume that $j\in\{1,2,\dots,k\}$,  $h_j$ verifies 
(\ref{eq:1}), and one of the following
assumptions is satisfied 
\begin{align}
\label{eq:15bi}
\mu_1(h_j)\geq -\Big(\frac{N-2}2\Big)^2+1 \quad&\text{and}\quad
\sum_{\substack{i=1\\i\not=j}}^k\frac{h_i\big(\frac{a_j-a_i}
  {|a_j-a_i|}\big)}{|a_j-a_i|^2}>0,\\
\label{eq:16bi}
-\Big(\frac{N-2}2\Big)^2<\mu_1(h_j)< -\Big(\frac{N-2}2\Big)^2+1
\quad&\text{and}\quad\sum_{\substack{i=1\\i\not=j}}^k
\int_{\R^N} \frac{h_i\big(\frac{x}{|x|}\big)\Big[\psi_1^{h_j}
 \big(\frac{x+a_i-a_j}{|x+a_i-a_j|}\big)\Big]^2}
   {|x|^2|x+a_i-a_j|^{2(\sigma_{h_j}+N-2)}}>0.
\end{align}
Then $S(h_1,\dots,h_k)<S(h_j)$.
\end{Proposition}
\begin{pf}
Since $h_j$ satisfies (\ref{eq:1}), by Theorem \ref{l:t96} there exists  
 $\phi_{h_j}\in\Di$, $\phi_{h_j}\geq0$, $\phi_{h_j}\not\equiv 0$, 
attaining $S(h_j)$, i.e. satisfying (\ref{eq:17bi}--\ref{eq:18bi}) with $h=h_j$.
Let us set $z_\mu(x)=\phi_\mu^{h_j}(x-a_j)$.
There holds
 \begin{align*}
   & S(h_1,\dots,h_k)\leq \frac {\displaystyle{\int_{\R^N} |\nabla
       z_\mu(x)|^2\,dx- \int_{\R^N}\frac{h_j\big(\frac{x-a_j}
         {|x-a_j|}\big)}{|x-a_j|^2}\,z_\mu^2(x)\,dx - \sum_{i\neq j}
       \int_{\R^N}\frac{h_i\big(\frac{x-a_i}
         {|x-a_i|}\big)}{|x-a_i|^2}\,z_\mu^2(x)\,dx}}
   {\|z_\mu\|_{L^{2^*}(\R^N)}^{2}}
   \\[10pt]
   &= \frac {\displaystyle{\int_{\R^N} |\nabla \phi_{h_j}(x)|^2\,dx-
       \int_{\R^N}\frac{h_j\big(\frac{x}
         {|x|}\big)}{|x|^2}\,\phi_{h_j}^2(x)\,dx - \sum_{i\neq j}
       \int_{\R^N}\frac{h_i\big(\frac{x-(a_i-a_j)}
         {|x-(a_i-a_j)|}\big)}{|x-(a_i-a_j)|^2}\,(\phi_\mu^{h_j}(x))^2\,dx}}
   {\|\phi_{h_j}\|_{L^{2^*}(\R^N)}^{2}}.
 \end{align*}
 From above and Lemmas \ref{l:interest1},
\ref{l:interest2}, and  \ref{l:interest3}, we deduce the following estimate
 \begin{align}\label{eq:45}
   & S(h_1,\dots,h_k)\leq 
S(h_j)-\|\phi_{h_j}\|_{L^{2^*}(\R^N)}^{-2}\times\\[8pt]
   &\notag\ \times\begin{cases}
     \mu^2\big(\int_{\R^N}\phi_{h_j}^2(x)\big)\Bigg(
     {\displaystyle{\sum_{i\neq j}}} \frac{h_i\big(\frac{a_j-a_i}
       {|a_j-a_i|}\big)}{|a_j-a_i|^2}+o(1)\Bigg)
     &\text{if }\mu_1(h_j)>-\big(\frac{N-2}2\big)^2+1\\[18pt]
     \mu^2\big(\int_{|x|<\frac1\mu}\phi_{h_j}^2(x)\big)\Bigg(
     {\displaystyle{\sum_{i\neq j}}} \frac{h_i\big(\frac{a_j-a_i}
       {|a_j-a_i|}\big)}{|a_j-a_i|^2}+o(1)\Bigg)
     &\text{if }\mu_1(h_j)=-\big(\frac{N-2}2\big)^2+1\\[18pt]
     \mu^{2\sigma_{h_j}+N-2}(c_{\infty}^{h_j})^2\Bigg(
     {\displaystyle{\sum_{i\not=j} \int_{\R^N}}}
     \frac{h_i(\frac{x}{|x|})\big[\psi_1^{h_j}
       \big(\frac{x+a_i-a_j}{|x+a_i-a_j|}\big)\big]^2}
     {|x|^2|x+a_i-a_j|^{2(\sigma_{h_j}+N-2)}}+o(1)\Bigg) &\text{if
     }\mu_1(h_j)<-\big(\frac{N-2}2\big)^2+1
\end{cases}
\end{align}
as $\mu\to0^+$.  Taking $\mu$ small enough in (\ref{eq:45}), from
(\ref{eq:15bi}--\ref{eq:16bi}) we obtain that $S(h_1,\dots,h_k)<S(h_j)$.
\end{pf}

\begin{remark}\label{rem:on}
  For $-\big(\frac{N-2}2\big)^2<\mu_1(h_j)<
  -\big(\frac{N-2}2\big)^2+1$, assumption (\ref{eq:16bi}) can be
  rewritten as
$$
\int_{\R^N}\Bigg( \sum_{i\not=j} \frac{(h_i\circ
  A_{ij})\big(\frac{x}{|x|}\big)}{|a_i-a_j|^{2(\sigma_{h_j}+N-2)}}
\Big[(\psi_1^{h_j}\circ A_{ij})
\Big(\frac{x+e_1}{|x+e_1|}\Big)\Big]^2\Bigg)
\frac{dx}{|x|^2|x+e_1|^{2(\sigma_{h_j}+N-2)}}>0,
$$
where $A_{ij}\in O(N)$ are such that $A_{ij}
e_1=\frac{a_i-a_j}{|a_i-a_j|}$.
\end{remark}

\section{The Palais-Smale condition and proof of Theorem
  \ref{t:ach}}\label{sec:pala-smale-cond}

\noindent 
If $u\in \Di$, $u>0$ a.e. in $\R^N$, is a critical points of the 
functional $J:\Di\to\R$,
\begin{equation}\label{eq:energy} 
J(v)=
\frac12\int_{\R^N} |\n v|^2dx-\frac12\sum_{i=1}^k
\int_{\R^N}
\dfrac{h_i\big(\frac{x-a_i}
{|x-a_i|}\big)}{|x-a_i|^2}v^2(x)\,dx-
\frac{S(h_1,h_2,\dots, h_k)}{2^*}\int_{\R^N} |v|^{2^*}dx,
\end{equation}
then $w=S(h_1,h_2,\dots, h_k)^{1/(2^*-2)}u$ is a solution to
equation~(\ref{eq:14}) (weakly in $\Di$ and classically in
$\R^N\setminus \{a_1,\dots,a_k\}$).  From now on, for any $u\in \Di$,
$J'(u)\in (\Di)^{\star}$ will denote the Fr\'echet derivative of $J$ at $u$
and $\langle\cdot,\cdot\rangle$ will stay for  the duality product
between $\Di$ and its dual space $(\Di)^\star$.

The Concentration-Compactness
analysis of the behavior of Palais-Smale sequences provides the
following local compactness result.
\begin{Theorem}\label{th:ps}
Let (\ref{eq:posde}) hold and  $\{u_n\}_{n\in\N}\subset \Di$
be a Palais-Smale sequence for
$J$, namely 
$$ 
\lim_{n\to\infty}J(u_n)=c<\infty\text{ in }\R\quad\text{and}\quad
\lim_{n\to\infty}J'(u_n)=0\text{ in the dual space }(\Di)^{\star}.
$$ 
If 
\begin{equation}\label{eq:19}
c<\frac{1}{N}S(h_1,h_2,\dots,h_k)^{1-\frac{N}{2}}\left(\min\bigg\{S,
S(h_1), \dots, S(h_k), S\Big(\sum\nolimits_{j=1}^k h_j\Big)\bigg\}\right)^{\!\!N/2},
\end{equation}
then $\{u_n\}_{n\in\N}$ admits a  subsequence strongly  converging in $\Di$.
\end{Theorem}
\begin{pf}
Let $\{u_n\}_{n\in\N}$ be a Palais-Smale sequence for $J$ at level $c$, then from
(\ref{eq:posde}) there exists some positive constant $c_1$ such that
\begin{align*}
c_1\|u_n\|_{\Di}^2&\leq 
{\mathcal Q}(u_n)
=NJ(u_n)-\frac{N-2}2
\langle J'(u_n),u_n\rangle =N c+o(\|u_n\|_{\Di})+o(1)
\end{align*}
as $n\to+\infty$, hence $\{u_n\}_{n\in\N}$ is a bounded sequence in $\Di$. Then
there exists $u_0\in\Di$ such that, up to
a subsequence still denoted as $\{u_n\}_{n\in\N}$, 
$u_n\rightharpoonup u_0$ weakly in $\Di$,
$u_n\to u_0$  a.e. in $\R^N$, and
$u_n\to u_0$ in  $L^\alpha _{\rm loc}(\R^N)$ for any  $\a \in [1,2^*)$.
The {\sl Concentration Compactness Principle} by P. L. Lions, (see
\cite{PL1} and \cite{PL2}), ensures  that, for an at most
countable set ${\mathcal J}$, some points $x_j\in\R^N\setminus\{a_1,\dots,a_k\}$,
some real numbers $\mu_{x_j},\nu_{x_j}$, $j\in{\mathcal
  J}$, and $\mu_{a_i},\nu_{a_i},\g_i$, $i=1,\dots,k$,  the
following convergences hold in the sense 
of measures up to a subsequence
\begin{align}
\label{eq:4}&|\nabla u_n|^2\rightharpoonup d\mu\ge |\n
u_0|^2+\sum_{i=1}^k \mu_{a_i}\delta_{a_i}+\sum_{j\in {\mathcal
    J} }\mu_{x_j}\delta_{x_j},\\
\label{eq:5}&|u_n|^{2^*}\rightharpoonup d\nu= |u_0|^{2^*}+\sum_{i=1}^k
\nu_{a_i}\delta_{a_i}+ \sum_{j\in {\mathcal
    J}}\nu_{x_j}\delta_{x_j},\\
\label{eq:6}& h_i\Big(\frac{x-a_i}
{|x-a_i|}\Big)\dfrac{u_n^2}{|x-a_i|^2}\rightharpoonup d \g_{a_i}=
h_i\Big(\frac{x-a_i}
{|x-a_i|}\Big)\dfrac{u_0^2}{|x-a_i|^2}+\g_i\delta_{a_i},\quad\text{for any}\quad
i=1,\dots,k. 
\end{align}
Notice that we can choose $\mu_{a_i},\mu_{x_j}$ such that
$\mu_{a_i}=d\mu(\{a_i\})$, $\mu_{x_j}=d\mu(\{x_j\})$. 
From Sobolev's inequality it follows that 
\begin{equation}\label{eq:8}
S\nu_{x_j}^{\frac2{2^*}}\le\mu_{x_j} \mbox{ for all } j\in{\mathcal
  J}\quad\text{and}\quad S\nu_{a_i}^{\frac2{2^*}}\le\mu_{a_i} \mbox{ for all }
i=1,\dots,k.
\end{equation}
The  concentration at infinity of  the sequence can valuated by  the following  quantities 
$$
\nu_{\infty}=\lim_{R\to \infty}\limsup_{n\to
\infty}\int_{|x|>R}|u_n(x)|^{2^*}dx,\quad \mu_{\infty}=\lim_{R\to
\infty}\limsup_{n\to \infty}\int_{|x|>R}|\n u_n(x)|^{2}dx$$ and $$
\g_{\infty}=\lim_{R\to \infty}\limsup_{n\to
\infty}\int_{|x|>R}\bigg(\sum_{i=1}^k h_i
\big({\textstyle{\frac{x}{|x|}}}\big)\bigg)\frac{ u_n^2(x)}{|x|^2}\,dx.
$$ 
Testing $J'(u_n)$ with $u_n\phi_j^\e$, for some smooth  cut-off
function $\phi_j^\e$ centered at $x_j$ and supported in $B(x_j,\e)$, and letting 
$n\to\infty$ and $\e\to 0$, we obtain that 
$\mu_{x_j} \leq S(h_1,h_2,\dots,h_k)\nu_{x_j}$, which, together with 
\eqref{eq:8}, implies that 
\begin{equation}\label{eq:7}
\mathcal J\quad\text{is finite and for }j\in {\mathcal J}\text{ 
either }\nu_{x_j}=0\text{ or }\nu_{x_j}\geq \bigg(\frac{S}{S(h_1,h_2,\dots,
h_k)}\bigg)^{\!\!N/2}.
\end{equation}
To analyze concentration at singularities, for each $i=1,2,\dots,k$ we
consider a smooth cut-off function $\psi_i^\e$ satisfying $0\leq
\psi_i^\e(x)\leq 1$,
$$
\quad \psi_i^\e(x)= 1\quad\mbox{if }|x-a_i|\leq \frac{\e}{2},\quad 
\psi_i^\e(x)= 0\quad\mbox{if }|x-a_i|\geq \e,\quad\text{and}\quad |\nabla
\psi_i^\e(x) |\leq \frac{4}{\e} \quad\mbox{for all }x\in\R^N. \quad
$$
From \eqref{eq:2} it follows that
that
\begin{equation*}
\frac{\int_{\R^N} |\n (u_n\psi_i^\e)|^2\,dx -\int_{\R^N}h_i\big(\frac{x-a_i}
{|x-a_i|}\big)
\frac{|\psi_i^\e|^2u^2_n}{|x-a_i|^2}\,dx}{\Big(\int_{\R^N}|\psi_i^\e
u_n|^{2^*}dx\Big)^{2/2^*}}\geq
S(h_i)
\end{equation*}
and hence
\begin{align}\label{eq:32}
  & \int_{\R^N} |\psi_i^\e| ^2|\n u_n|^2dx +\int_{\R^N} u_n^2|\n
  \psi_i^\e|^2\,dx+2\int_{\R^N}
  u_n\psi_i^\e \n u_n\cdot \n \psi_i^\e dx  \\
  \notag&\hskip2cm \geq\int_{\R^N}h_i\Big(\frac{x-a_i} {|x-a_i|}\Big)
  \frac{|\psi_i^\e|^2u^2_n}{|x-a_i|^2}\,dx+
  S(h_i)\bigg(\int_{\R^N}|\psi_i^\e u_n|^{2^*}dx\bigg)^{\!\!2/2^*}.
\end{align}
It is easy to verify that 
\begin{equation*}
\lim_{\e\to 0}\limsup_{n\to \infty}\bigg[\int_{\R^N}
u_n^2|\n \psi_i^\e|^2\,dx+2\int_{\R^N} u_n \psi_i^\e \n u_n\cdot\n \psi_i^\e
\,dx\bigg]=0,
\end{equation*}
then from (\ref{eq:32})  and (\ref{eq:4}--\ref{eq:6}) we deduce that
\begin{equation}\label{eq:10}
\mu_{a_i}\geq \g_i+
S(h_i)\nu_{a_i}^{2/2^*}.
\end{equation}
Testing $J'(u_n)$ with $u_n\psi_i^\e$ and letting $n\to+\infty$ and
$\e\to0$, we obtain that
\begin{equation}\label{eq:11}
\mu_{a_i}-\g_{i}\leq S(h_1,h_2,\dots,h_k)\nu_{a_i}. 
\end{equation}
From (\ref{eq:10}) and (\ref{eq:11}) we 
conclude that, for each $i=1,2,\dots,k$,
\begin{equation}\label{eq:9}
\text{either}\quad\nu_{a_i}=0\quad\text{or}\quad \nu_{a_i}\geq
\bigg(\frac{S(h_i)}{S(h_1,h_2,\dots,h_k)}\bigg)^{\!\!N/2}.
\end{equation}
To  study the possibility of concentration at  $\infty$, we consider 
 a regular cut-off function  $\psi_R$ such that 
$$
0\le \psi_R(x)\le 1\text{ for all }x\in\R^N,\ 
\psi_R(x)=
\left\{
\begin{array}{l}
1,\,\mbox{ if }|x|>2R, \\ 0,\,\mbox{ if }|x|<R,
\end{array}
\right. \text{ and }|\nabla \psi_R(x) |\leq \frac2R\text{ for
  all }x\in \R^N .
$$
From \eqref{eq:2} we obtain
that
\begin{equation*}
\frac{\int_{\R^N} |\n (u_n\psi_R)|^2dx -\int_{\R^N}\Big(\sum_{i=1}^k
h_i\big(\frac{x}
{|x|}\big)\Big)
\frac{\psi_R^2u^2_n}{|x|^2}dx}{\Big(\int_{\R^N}|\psi_R
u_n|^{2^*}dx\Big)^{2/2^*}}\geq
S\Big({\textstyle{\sum_{i=1}^k h_i}}\Big)
\end{equation*}
and, consequently,
\begin{align}\label{eq:14b}
&\int_{\R^N} \psi_R ^2|\n u_n|^2dx +\int_{\R^N} u_n^2|\n \psi_R|^2dx+2\int_{\R^N}
u_n\psi_R \n u_n\cdot \n \psi_R\, dx \\
 \notag &\hskip2cm \geq  
\int_{\R^N}\Big({\textstyle{\sum_{i=1}^k
h_i\big(\frac{x}
{|x|}\big)}}\Big)
\frac{\psi_R^2u^2_n}{|x|^2}dx
+S\Big({\textstyle{\sum_{i=1}^k h_i}}\Big)
\bigg(\int_{\R^N}|\psi_R
u_n|^{2^*}dx\bigg)^{\!\!2/2^*}.
\end{align}
It is easy to verify that 
$$
\lim_{R\to \infty}\limsup_{n\to \infty}\bigg\{\int_{\R^N}
u_n^2|\n \psi_R|^2dx+2\int_{\R^N} u_n\psi_R \n u_n\cdot\n \psi_R\,
dx\bigg\}=0.
$$
Then from (\ref{eq:14b}) we infer
\begin{equation}\label{eq:15}
\mu_{\infty}-\g_{\infty}\geq
S\Big({\textstyle{\sum_{i=1}^kh_i}}\Big)\nu_{\infty}^{2/2^*}.
\end{equation}
Testing $J'(u_n)$ with $u_n\psi_R$ we obtain
\begin{align}\label{eq:21}
 0&=\lim_{n\to \infty }\langle
J'(u_n),u_n\psi_R\rangle  \\
\notag&=\lim_{n\to \infty}\bigg[\int_{\R^N} |\nabla
  u_n|^2\psi_R+\int_{\R^N} u_n \nabla u_n \cdot\nabla 
 \psi_R \\
\notag&\qquad\qquad- 
\sum_{i=1}^k\int_{\R^N}
h_i\Big(\frac{x-a_i}
{|x-a_i|}\Big)\frac{\psi_Ru^2_n}{|x-a_i|^2}dx
 -S(h_1,h_2,\dots,
h_k)\int_{\R^N} \psi_R
|u_n|^{2^*}\bigg].
\end{align}
If $|x|\geq R$ with $R$ sufficiently large, there holds
\begin{align*}
  \bigg| \frac{h_i\big(\frac{x-a_i} {|x-a_i|}\big)}{|x-a_i|^2}&-
  \frac{h_i\big(\frac{x} {|x|}\big)}{|x|^2}\bigg|\leq \left|
    \frac{h_i\big(\frac{x-a_i} {|x-a_i|}\big)}{|x-a_i|^2}-
    \frac{h_i\big(\frac{x-a_i} {|x-a_i|}\big)}{|x|^2}\right|+
  \frac1{|x|^2}\,\left| h_i\Big(\frac{x-a_i}
    {|x-a_i|}\Big)-h_i\Big(\frac{x}
    {|x|}\Big)\right|\\
  &\leq \|h_i\|_{L^{\infty}({\mathbb S}^{N-1})}\frac{|2a_i\cdot
    x-|a_i|^2|}{|x-a_i|^2|x|^2} +\frac{{\rm const}}{|x|^2} \,\left|
    \frac{x-a_i}
    {|x-a_i|}-\frac{x}{|x|}\right|\\
  &\leq \|h_i\|_{L^{\infty}({\mathbb
      S}^{N-1})}\frac{2|a_i||x|+|a_i|^2}{|x-a_i|^2|x|^2} +\frac{{\sqrt
      2\,\rm const}}{|x|^2} \,\sqrt{\frac{|x|\big(|x-a_i|-|x|\big)+a_i\cdot
    x}{|x-a_i||x|}}\leq \frac{{\rm const}}{|x|^{5/2}}.
\end{align*}
Since, by H\"older's inequality,
$$
\int_{\R^N}\frac{ u_n^2\psi_R}{|x|^{5/2}}\,dx \leq
\bigg(\int_{|x|>R}u_n^{2^*}\bigg)^{\!\!2/2^*} \bigg(\int_{|x|>R}|x|^{-\frac54N}\bigg)^{\!\!2/N} 
 =O(R^{-1/2})
$$
as $R\to+\infty$ uniformly with respect to $n$,
we deduce that 
$$
\sum_{i=1}^k\int_{\R^N} h_i\Big(\frac{x-a_i}
{|x-a_i|}\Big)\frac{\psi_Ru^2_n}{|x-a_i|^2}dx= \int_{\R^N}
\frac{\sum_{i=1}^k h_i\big(\frac{x}
  {|x|}\big)}{|x|^2}\psi_Ru^2_n\,dx+O(R^{-1/2})
$$
as $R\to+\infty$ uniformly with respect to $n$, hence
\begin{equation}\label{eq:16}
\lim_{R\to+\infty}\limsup_{n\to\infty}
\int_{\R^N}
\sum_{i=1}^kh_i\Big(\frac{x-a_i}
{|x-a_i|}\Big)\frac{\psi_Ru^2_n}{|x-a_i|^2}dx=
\g_{\infty}.
\end{equation}
Passing to lim-sup as $n\to\infty$ and limits as $R\to\infty$ in
(\ref{eq:21}) and using (\ref{eq:16}), we obtain  that
\begin{equation}\label{eq:17}
\mu_{\infty}-\g_{\infty}= S(h_1,h_2,\dots,h_k)\nu_{\infty}. 
\end{equation}
From (\ref{eq:15}) and (\ref{eq:17}) we 
conclude that 
\begin{equation}\label{eq:12bi}
\text{either}\quad\nu_{\infty}=0\quad\text{or}\quad \nu_{\infty}\geq
\bigg(\frac{S(\sum_{i=1}^k h_i)}{S(h_1,h_2,\dots,h_k)}\bigg)^{\!\!{N}/{2}} .
\end{equation}
As a conclusion we obtain
\begin{align}\label{eq:18} 
c & =  J(u_n)-\frac{1}{2}\langle J'(u_n),u_n\rangle
  +o(1) =  \frac{1}{N}S(h_1,h_2,\dots,h_k)\int_{\R^N}
|u_n|^{2*}dx +o(1)\\
& = \frac{S(h_1,h_2,\dots,h_k)}{N}\bigg\{\int_{\R^N}
|u_0|^{2*}dx+\sum_{i=1}^k\nu_{a_i}+\nu_{\infty}+\sum_{j\in {\mathcal
J}}\nu_{x_j}\bigg\}.\notag
\end{align}
From (\ref{eq:19}), (\ref{eq:18}), (\ref{eq:7}), (\ref{eq:9}),
and (\ref{eq:12bi}), we deduce that $\nu_{x_j}=0$ for any $j\in {\mathcal
  J}$, $\nu_{a_i}=0$ for any $i=1,\dots,k$, and $\nu_{\infty}=0$. Then
up to a subsequence $u_n\to u_0$ in $\Di$.
\end{pf}

The Palais-Smale condition recovered in Theorem \ref{th:ps} and the interaction estimates 
proved in Proposition \ref{cor:sub} are the key tools to prove Theorem \ref{t:ach}.

\medskip\noindent
\begin{pfn}{Theorem \ref{t:ach}}
  Let $\{u_n\}_n\subset\Di$ be a minimizing sequence for
  \eqref{eq:min}. From the homogeneity of the quotient, we can require
  without restriction that $\|u_n\|_{L^{2^*}(\R^N)}=1$, while from
  Ekeland's variational principle we can assume that the sequence satisfies
  the Palais-Smale property, i.e. for any $v\in\Di$
\begin{multline*}
\int_{\R^N}\n u_n(x)\cdot \n v(x)\,dx
-\sum_{i=1}^k
\int_{\R^N}
\dfrac{h_i\big(\frac{x-a_i}
{|x-a_i|}\big)u_n(x)}{|x-a_i|^2}v(x)\,dx
\\-
S(h_1,h_2,\dots, h_k)
\int_{\R^N}
|u_n(x)|^{2^*-2}u_n(x)v(x)\,dx=o\big(\|v\|_{\Di}\big).
\end{multline*}
Hence $J'(u_n)\to 0$ in $(\Di)^{\star}$ and 
$$
J(u_n)\to \Big(\frac12-\frac1{2^*}\Big) S(h_1,h_2,\dots, h_k)=\frac1N
S(h_1,h_2,\dots, h_k).
$$
From assumption \eqref{eq:25} and Proposition \ref{cor:sub}, we infer
that 
 \begin{equation}\label{eq:36}
S(h_1,h_2,\dots,h_k)<S(h_k).
\end{equation}
From assumptions
\eqref{eq:23} and (\ref{eq:26})  we have that
 \begin{equation}\label{eq:37}
S(h_k)\leq S(h_i)\quad\text{for all}\quad i=1,\dots,k-1,
\quad\text{and}\quad S(h_k)\leq S\Big({\textstyle{\sum_{i=1}^k}} h_i  \Big),
\end{equation}
while from assumption (\ref{eq:24}) and Theorem \ref{l:t96} there holds
\begin{equation}\label{eq:34}
 S(h_k)< S.
\end{equation}
Gathering  (\ref{eq:36}), (\ref{eq:37}), and (\ref{eq:34}), we finally
have
$$
S(h_1,h_2,\dots,h_k)<
\min\bigg\{S,
S(h_1), \dots, S(h_k), S\Big({\textstyle{\sum_{i=1}^k}} h_i  \Big)\bigg\}
$$
and hence 
$$
\frac1N
S(h_1,h_2,\dots,
h_k)<\frac{1}{N}S(h_1,h_2,\dots,h_k)^{1-\frac{N}{2}}\left(\min\bigg\{S, 
S(h_1), \dots, S(h_k),
S\Big(\sum\nolimits_{i=1}^kh_i\Big)\bigg\}\right)^{\!\!N/2} .
$$
From Theorem \ref{th:ps}  we 
deduce that $\{u_n\}_{n\in\N}$ has a  subsequence strongly
converging to some $u_0\in\Di$, which satisfies
$J(u_0)=\frac1NS(h_1,h_2,\dots,h_k)$.
In particular $u_0$ achieves the infimum in \eqref{eq:min}. Since
$J(u_0)=J(|u_0|)$, we have that also $|u_0|$ is a minimizer in \eqref{eq:min}
and then $v_0=S(h_1,h_2,\dots,h_k)^{1/(2^*-2)}|u_0|$ is a
nonnegative solution to equation \eqref{eq:14}. The maximum principle
in $\R^N\setminus\{a_1,\dots,a_k\}$ implies that $v_0>0$ in $\R^N\setminus\{a_1,\dots,a_k\}$.
\end{pfn}

Let us now consider the case of singularities generated by electric
dipoles.  In order to prove Corollary \ref{c:ach}, we first need to
establish the following monotonicity property of ground state levels
with respect to the dipole magnitudes.
\begin{Lemma}\label{l:mono}
  If $\Lambda_N^{-1}>\lambda_1\geq\lambda_2>0$, ${\mathbf d}_1,{\mathbf d}_2\in\R^N$
  with $|{\mathbf d}_1|=|{\mathbf d}_2|=1$, and
  $h_i(\theta)=\lambda_i\theta\cdot {\mathbf d}_i$ for $i=1,2$, then
$S(h_2)\geq S(h_1)$.
\end{Lemma}
\begin{pf}
We first notice that, by rotation invariance, for any $\lambda>0$,
$S(\lambda \theta \cdot {\mathbf d})$ 
does not depend on the unit vector ${\mathbf d}$, hence
$S(h_2)=S(\tilde{h}_2)$ where 
$\tilde h_2(\theta)=\lambda_2\theta\cdot {\mathbf d}_1$.

From Theorem \ref{l:t96},  there exists $w\in\Di\setminus\{0\}$ such that
\begin{equation}\label{eq:33}
\frac{\int_{\R^N}\big[
    |\nabla w(x)|^2-\frac{\lambda_2 x\cdot  {\mathbf d}_1}{|x|^3}w^2(x)\big]\,dx}
  {\big(\int_{\R^N}|w(x)|^{2^*}dx\big)^{2/2^*}}=S(\tilde h_2).
\end{equation}
We claim that the  quotient at the left hand side decreases after passing to polarization  with
respect to the half-space $H_{{\mathbf d_1}}:=\{x\in\R^N:\,x\cdot
{\mathbf d_1}\geq0\}$. We denote as $\sigma_{\mathbf
      d_1}:\R^N\to\R^N$ the reflection with respect to the boundary of
$H_{{\mathbf d_1}}$, i.e. $\sigma_{{\mathbf
      d_1}}(x)=x-2(x\cdot{\mathbf d_1}){\mathbf d_1}$. The polarization of any measurable
nonnegative function $u$ with respect to $H_{{\mathbf d_1}}$ is defined as
$$
 u_{{\mathbf d_1}}(x):=\begin{cases}
\max\{u(x),u(\sigma_{{\mathbf d_1}}(x))\},&\text{if }x\in H_{{\mathbf d_1}},\\
\min\{u(x),u(\sigma_{{\mathbf d_1}}(x))\},&\text{if }x\in \R^N\setminus
H_{{\mathbf d}_1}.
\end{cases}
$$
From well known properties of polarization, there holds
\begin{equation}\label{eq:35}
\|\nabla |w|_{{\mathbf d_1}}\|_{L^2(\R^N)}=\|\nabla w\|_{L^2(\R^N)}
\quad\text{and}\quad
\| |w|_{{\mathbf d_1}}\|_{L^{2^*}(\R^N)}=\| w\|_{L^{2^*}(\R^N)},
\end{equation}
see
\cite[Propositions 22.2 and 22.5]{Willem}. 
Moreover 
\begin{align}\label{eq:38}
&\int_{\R^N}\frac{x\cdot{\mathbf d_1} }
{|x|^3}\,\big(
|w|_{{\mathbf d_1}}^2-w^2\big)\,dx\\
&\notag\qquad=\int_{ H_{{\mathbf d_1}}}\!\!\frac{x\cdot{\mathbf d_1} }
{|x|^3}\,\big(
|w|_{{\mathbf d_1}}^2-|w|^2\big)\,dx+
\int_{\R^N\setminus H_{{\mathbf d_1}}}\!\!\frac{x\cdot{\mathbf d_1} }
{|x|^3}\,\big(
|w|_{{\mathbf d_1}}^2-|w|^2\big)\,dx\geq0
\end{align}
and, through the change of variables $x=\sigma_{{\mathbf d_1}}(y)$,
\begin{align}\label{eq:39}
  \int_{\R^N}\frac{x\cdot{\mathbf d_1} } {|x|^3}\,|w|_{{\mathbf
      d_1}}^2(x)\,dx&=\int_{ H_{{\mathbf d_1}}}
  \!\!\frac{x\cdot{\mathbf d_1} } {|x|^3}\,|w|_{{\mathbf
      d_1}}^2(x)\,dx+ \int_{\R^N\setminus H_{{\mathbf d_1}}}\!\!
  \frac{x\cdot{\mathbf d_1} }
  {|x|^3}\,|w|_{{\mathbf d_1}}^2(x)\,dx\\
  \notag&=\int_{ H_{{\mathbf d_1}}} \!\!\frac{x\cdot{\mathbf d_1} }
  {|x|^3}\,|w|_{{\mathbf d_1}}^2(x)\,dx- \int_{H_{{\mathbf d_1}}}\!\!
  \frac{y\cdot{\mathbf d_1} }
  {|y|^3}\,|w|_{{\mathbf d_1}}^2(\sigma_{{\mathbf d_1}}(y))\,dy\\
  \notag&=\int_{ H_{{\mathbf d_1}}} \!\!\frac{x\cdot{\mathbf d_1} }
  {|x|^3}\,\big( |w|_{{\mathbf d_1}}^2(x)- |w|_{{\mathbf
      d_1}}^2(\sigma_{{\mathbf d_1}}(x))\big)\,dx\geq0.
\end{align}
From (\ref{eq:33}--\ref{eq:39}), we obtain that
\begin{align*}
S(h_2)=S(\tilde h_2)&\geq \frac{\int_{\R^N}\big[ |\nabla |w|_{{\mathbf
      d_1}}(x)|^2-\frac{\lambda_2 x\cdot {\mathbf
      d}_1}{|x|^3}|w|_{{\mathbf d_1}}^2(x)\big]\,dx}
{\big(\int_{\R^N}||w|_{{\mathbf d_1}}(x)|^{2^*}dx\big)^{2/2^*}}\\
&\geq \frac{\int_{\R^N}\big[ |\nabla |w|_{{\mathbf
      d_1}}(x)|^2-\frac{\lambda_1 x\cdot {\mathbf
      d}_1}{|x|^3}|w|_{{\mathbf d_1}}^2(x)\big]\,dx}
{\big(\int_{\R^N}||w|_{{\mathbf d_1}}(x)|^{2^*}dx\big)^{2/2^*}}
\geq S(h_1),
\end{align*}
thus proving the stated inequality.
\end{pf}

\medskip\noindent
\begin{pfn}{Corollary \ref{c:ach}}
Theorem \ref{t:ach} applies with $h_i(\theta)=\lambda_i\theta\cdot {\mathbf d}_i$. 
Indeed, (\ref{eq:23}) follows from Lemma \ref{l:mono}, (\ref{eq:25}) comes from  
(\ref{eq:27}) and (\ref{eq:26}) from (\ref{eq:29}) and Lemma \ref{l:mono}. 
\end{pfn}

\section{The problem on bounded domains}\label{sec:probl-bound-doma}

\noindent In this
section we discuss the existence of ground state solutions to  \eqref{eq:bound}
by analyzing  the associated  minimization problem
 (\ref{eq:minbound}) on a bounded smooth domain $\Omega\subset\R^N$, $N\geq 3$,
containing points $a_1,\dots,a_k$. 
The corresponding functional is given by
\begin{equation}\label{eq:energybounded} 
J_{\Omega}(v)=
\frac12\int_{\Omega} |\n u|^2dx
-\frac12\sum_{i=1}^k
\int_{\R^N}
\dfrac{h_i\big(\frac{x-a_i}
{|x-a_i|}\big)}{|x-a_i|^2}v^2(x)\,dx-
\frac{S(h_1,h_2,\dots, h_k)}{2^*}
\int_{\Omega} |v|^{2^*}dx.
\end{equation}
By boundedness of the domain, minimizing sequences of  (\ref{eq:minbound}) 
cannot lose  mass at infinity. Hence, arguing as in Theorem
\ref{th:ps}, the following  local Palais-Smale condition can be obtained.
\begin{Theorem}\label{th:psbound}
Assume that \eqref{eq:pdbd} holds. Let $\{u_n\}_{n\in\N}\subset H^1_0(\Omega)$
be a Palais-Smale sequence for
$J_{\Omega}$, namely $\lim_{n\to\infty}J_{\Omega}(u_n)=c$ in $\R$ and 
$\lim_{n\to\infty}J_{\Omega}'(u_n)=0$ in the dual space $(H_0^1(\Omega))^{\star}$.  
If 
\begin{equation*}
c<c_{\Omega}^*=\frac{1}{N}S_{\Omega}(h_1,h_2,\dots,h_k)^{1-\frac{N}{2}}\min\bigg\{S,
S(h_1), \dots, S(h_k)\bigg\}^{\!\!N/2},
\end{equation*}
then $\{u_n\}_{n\in\N}$ has a converging subsequence.
\end{Theorem}

In a bounded domain, the comparison between ground state levels of
dipole-type and multi-dipole type problems is more delicate and
requires an analysis of the concentration behavior of cutted-off test
functions. To this aim we need, besides the asymptotic behavior of
functions $\phi_h$ at infinity, also the behavior of their
gradient, which we are going to deduce from  Green's representation formula
and the following  property of
differentiability of {\em{Newtonian potentials}}.

\begin{Lemma}\label{l:green}
  Let $\Omega\subset\R^N$ be a bounded smooth domain, $g\in
  L^p(\Omega)$, for every $p\in [1,2)$, and let $u$ be the Newtonian
  potential of $g$, i.e.
$$
u(x)=\frac1{N(2-N)\omega_N}\int_{\Omega}\frac{g(y)}{|x-y|^{N-2}}\,dy.
$$
Then $u\in W^{1,q}(\R^N)$ for all
$q\in\big(\frac{N}{N-2},\frac{2N}{N-2}\big)$ and  the weak derivatives of
$u$ are given by 
$$
\frac{\partial u}{\partial
  x_i}(x)=\frac1{N\omega_N}
\int_{\Omega}\frac{g(y)(x_i-y_i)}{|x-y|^{N}}\,dy,\quad i=1,\dots,N. 
$$
\end{Lemma}

\begin{pf}
The proof can be obtained by approximation from \cite[Lemma 4.1,
p. 54]{GT} using  the $L^p$ inequalities  for singular Riesz potentials 
proved in \cite[Theorem 1, p. 119]{stein}. We refer to \cite[Lemma A.1]{FMT}
for a detailed proof in the case $g\in L^2(\Omega)$, which can be followed
step by step yielding Lemma~\ref{l:green}.
\end{pf}

From the above lemma and  Green's representation formula we derive the following estimate 
on the behavior of solutions $\phi_h$ as $|x|\to+\infty$.
 
\begin{Lemma}\label{lem:stimagrad}
For  $h\in C^1({\mathbb S}^N)$ verifying
$\mu_1(h)\geq-\big(\frac{N-2}2\big)^{\!2}+1$
 and (\ref{eq:1}), let $\phi_h\in\Di$,
$\phi_h\geq0$, $\phi_h\not\equiv 0$, be as in (\ref{eq:17bi}--\ref{eq:18bi}).
Then, for every $\e>0$,
\begin{equation}\label{eq:stimagradiente}
|\nabla \phi_h(x)|=
\begin{cases}
O\big(|x|^{-\sigma_h-N+1}\big),&\text{if }\mu_1(h)<N-1,\\
O\big(|x|^{-N+\e}\big),&\text{if }\mu_1(h)\geq N-1,
\end{cases}
\quad\text{as }
|x|\to+\infty.
\end{equation}
\end{Lemma}
\begin{pf}
  Let $w_h(x):=|x|^{-(N-2)}\phi_h(x/|x|^2)$ be the Kelvin transform
  of $\phi_h$.  Then $w_h$ solves 
$$
-\Delta   w_h= g \quad\text{in }\R^N,
$$
where 
$$
g(x)=\frac{h(x/|x|)}{|x|^2}w_h(x)+w_h^{2^*-1}(x).
$$
Moreover, a direct calculation yields the following relation between the gradients 
of $\phi_h$ and of its Kelvin transform
\begin{equation}\label{eq:44}
\nabla \phi_h(x)=|x|^{-N}\nabla w_h\Big(\frac x{|x|^2}\Big)-2x|x|^{-N-2}x\cdot\nabla 
w_h\Big(\frac x{|x|^2}\Big)-(N-2)|x|^{-N}w_h\Big(\frac x{|x|^2}\Big)x.
\end{equation}
From (\ref{eq:7bi}), $ w_h(x)=O(|x|^{\sigma_h})$  as $x\to 0$, 
hence  $g(x)=O(|x|^{\sigma_h-2})$ as $x\to 0$.
Therefore, from  $\mu_1(h)\geq-\big(\frac{N-2}2\big)^{\!2}+1$, 
it follows that $g\in L^p(B(0,1))$ for every $p\in [1,2)$.

Green's representation formula yields
\begin{align}\label{eq:43}
w_h(x)=&\frac{1}{N(N-2)\omega_N}\bigg[\int_{B(0,1)}
\frac{g(y)}{|x-y|^{N-2}}\,dy
+\int_{\partial
  B(0,1)}\frac{1}{|x-y|^{N-2}}\,\frac{\partial w_h}{\partial
  \nu}\,dS(y)\bigg]\\
\notag&+\frac{1}{N\omega_N}\int_{\partial
  B(0,1)}\frac{w_h(y)}{|x-y|^{N}}\, (y-x)\cdot\nu(y)\,dS(y),
\quad x\in B(0,1),
\end{align}
where $\omega_N$ denotes the volume of the unit ball in $\R^N$, 
$\nu$ is the unit outward normal to $\partial
  B(0,1)$, and $dS$ indicates the $(N-1)$-dimensional area
  element in $\partial   B(0,1)$.
It is easy to verify that the functions
\begin{align*}
x\mapsto \int_{\partial
  B(0,1)}\frac{1}{|x-y|^{N-2}}\,\frac{\partial w_h}{\partial
  \nu}\,dS(y),\qquad
 x\mapsto \int_{\partial
  B(0,1)}\frac{w_h(y)}{|x-y|^{N}}\, (y-x)\cdot\nu(y)\,dS(y),
\end{align*}
are of class $C^1( B(0,1))$. From Lemma \ref{l:green}, we have that
$$
\n\bigg(
\frac{1}{N(N-2)\omega_N}\int_{B(0,1)} 
\frac{g(y)}{|x-y|^{N-2}}\,dy
\bigg)=-\frac{1}{N\omega_N}\int_{B(0,1)}
\frac{x-y}{|x-y|^N}\,g(y)\,dy,
$$
and hence
\begin{align}\label{eq:40}
  \bigg|\n\bigg( \frac{1}{N(N-2)\omega_N}\int_{B(0,1)}
  \frac{g(y)}{|x-y|^{N-2}}\,dy \bigg)\bigg|\leq {\rm
    const\,}\int_{B(0,1)}\frac{|y|^{\sigma_h-2}}{|x-y|^{N-1}} \,dy.
\end{align}
If $\mu_1(h)<N-1$, i.e. $\sigma_h<1$, then 
$$
\bigg|\n\bigg(
\frac{1}{N(N-2)\omega_N}\int_{B(0,1)} 
\frac{g(y)}{|x-y|^{N-2}}\,dy 
\bigg)\bigg|\leq {\rm const\,}f(x),
$$
where $$
f(x)=\int_{\R^N}\frac{|y|^{\sigma_h-2}}{|x-y|^{N-1}}
\,dy.
$$
An easy scaling argument shows that $f(\alpha x)=\alpha^{\sigma_h-1}f(x)$
for all $\alpha>0$, hence $f(x)=|x|^{\sigma_h-1}f(e_1)$, where
$e_1=(1,0,\dots,0)\in\R^N$. Then, if $\mu_1(h)<N-1$,
\begin{align}\label{eq:41}
 \bigg|\n\bigg(
\frac{1}{N(N-2)\omega_N}\int_{B(0,1)} 
\frac{g(y)}{|x-y|^{N-2}}\,dy 
\bigg)\bigg|\leq {\rm const\,} |x|^{\sigma_h-1}.
\end{align}
If $\mu_1(h)\geq N-1$, i.e. $\sigma_h\geq1$,  we fix $0<\e<N-1$
and notice that, from (\ref{eq:40}),  
$$
\bigg|\n\bigg(
\frac{1}{N(N-2)\omega_N}\int_{B(0,1)} 
\frac{g(y)}{|x-y|^{N-2}}\,dy 
\bigg)\bigg|\leq {\rm const\,} k_{\e}(x),
$$
where 
$$
k_{\e}(x)=
\int_{\R^N} \frac1{|y|^{1+\e}|y-x|^{N-1}}\,dy.
$$
An easy scaling argument shows that $k_{\e}(\alpha x)=\alpha^{-\e}k_{\e}(x)$
for all $\alpha>0$, hence $k_{\e}(x)=|x|^{-\e}k_{\e}(e_1)$. Then, if $\mu_1(h)\geq N-1$
\begin{align}\label{eq:42}
 \bigg|\n\bigg(
\frac{1}{N(N-2)\omega_N}\int_{B(0,1)} 
\frac{g(y)}{|x-y|^{N-2}}\,dy 
\bigg)\bigg|\leq C(\e) |x|^{-\e},
\end{align}
for some positive constant $C(\e)$ depending on $\e$ (and also on $N$,
$h$, and $w_h$). Representation (\ref{eq:43}), regularity of the
boundary terms, and estimates (\ref{eq:41}--\ref{eq:42}) yield
\begin{equation}\label{eq:39bi}
\n w_h(x)
=
\begin{cases}
O\big(|x|^{\sigma_h-1}\big),&\text{if }\mu_1(h)<N-1,\\[5pt]
O\big(|x|^{-\e}\big),&\text{if }\mu_1(h)\geq N-1,
\end{cases}
\qquad \text{ as }x\to 0.
\end{equation}
Estimate (\ref{eq:stimagradiente}) follows then from (\ref{eq:39bi}) and
(\ref{eq:44}).  
\end{pf}

\begin{Lemma}\label{lem:sub}
Let $j\in\{1,2,\dots,k\}$. There holds
\begin{align}\label{eq:47}
&S_{\Omega}(h_1,\dots,h_k)\leq 
S(h_j)+O(\mu^{2\sigma_{h_j}+N-2})\\
\nonumber&-
\begin{cases}
  \mu^2\|\phi_{h_j}\|_{L^{2^*}(\R^N)}^{-2}\big(\int_{\R^N}\phi_{h_j}^2(x)\big)\Bigg(
  {\displaystyle{\sum_{i\neq j}}} \frac{h_i\big(\frac{a_j-a_i}
    {|a_j-a_i|}\big)}{|a_j-a_i|^2}+o(1)\Bigg),
  &\text{if }\mu_1(h_j)>-\big(\frac{N-2}2\big)^2+1,\\[18pt]
  \mu^2\|\phi_{h_j}\|_{L^{2^*}(\R^N)}^{-2}
  \big(\int_{|x|<\frac1\mu}\phi_{h_j}^2(x)\big)\Bigg(
  {\displaystyle{\sum_{i\neq j}}} \frac{h_i\big(\frac{a_j-a_i}
    {|a_j-a_i|}\big)}{|a_j-a_i|^2}+o(1)\Bigg), &\text{if
  }\mu_1(h_j)=-\big(\frac{N-2}2\big)^2+1,
\end{cases}
\end{align}
as $\mu\to 0^+$.
\end{Lemma}
\begin{pf}
  Let $\omega$ be an open set such that $\overline\omega\subset\Omega$
  and $a_j\in\omega$ and let $\psi\in C^{\infty}_c(\R^N)$ be a smooth
  cut-off function such that $0\leq \psi(x)\leq 1$, $\psi\equiv 0$ in
  $\R^N\setminus\Omega$, $\psi\equiv 1$ in $\omega$.  Then $\psi(x)
  \phi^{h_j}_{\mu}(x-a_j)\in H^1_0(\Omega)$.

Let $0<\e<\frac{N-2}2$. We claim that, as $\mu\to 0^+$, the following estimates hold:
\begin{align}
\label{eq:46}
\int_{\R^N}|\n
    (\psi(x)\phi^{h_j}_{\mu}(x-a_j))|^2\,dx&=
\int_{\R^N}|\n \phi_{h_j}
(x)|^2\,dx+ O(\mu^{2\sigma_{h_j}+N-2})+
O\big(\mu^{N-2\e}\big)\\[5pt]
\label{eq:48}
\int_{\R^N}
    \frac{h_j\big(\frac{x-a_j}
      {|x-a_j|}\big)}{|x-a_j|^2}|\psi(x)\phi^{h_j}_{\mu}(x-a_j)|^2
    \,dx&=
\int_{\R^N}
    \frac{h_j\big(\frac{x}
      {|x|}\big)}{|x|^2}|\phi_{h_j}(x)|^2 \,dx
+  O(\mu^{2\sigma_{h_j}+N-2})\\[5pt]
\label{eq:49}
\int_{\R^N}\dfrac{h_i\big(\frac{x-a_i}
      {|x-a_i|}\big)}{|x-a_i|^2}|\psi(x)\phi^{h_j}_{\mu}(x-a_j)|^2\,dx
&=
\int_{\R^N}\dfrac{h_i\big(\frac{x+a_j-a_i}
      {|x+a_j-a_i|}\big)}{|x+a_j-a_i|^2}|\phi_{\mu}^{h_j}(x)|^2\,dx
+  O(\mu^{2\sigma_{h_j}+N-2})\\[5pt]
\label{eq:51}
\bigg({\displaystyle{\int_{\R^N}}}
|\psi(x)\phi^{h_j}_{\mu}(x-a_j)|^{2^*}\,dx\bigg)^{\!\!2/2^*}&=
\bigg(\int_{\R^N} |\phi_{h_j}(x)|^{2^*}\,dx\bigg)^{\!\!2/2^*}+
O(\mu^{2\sigma_{h_j}+N-2}).
\end{align}
Let us prove (\ref{eq:46}). We have that
\begin{align}\label{eq:61}
\int_{\R^N}&|\n
(\psi(x)\phi^{h_j}_{\mu}(x-a_j))|^2\,dx=
\int_{\R^N}\psi^2(x)|\n \phi^{h_j}_{\mu}(x-a_j)|^2\,dx\\
\nonumber &+
\int_{\R^N}|\phi^{h_j}_{\mu}(x-a_j)|^2 |\n\psi(x)|^2\,dx+
2\int_{\R^N}\psi(x)\phi^{h_j}_{\mu}(x-a_j)\n\psi(x)\cdot \n
\phi^{h_j}_{\mu}(x-a_j)\, dx .
\end{align}
In view of \eqref{eq:stimagradiente} we have
\begin{align}\label{eq:67}
&\bigg|\int_{\R^N}\psi^2(x)|\n \phi^{h_j}_{\mu}(x-a_j)|^2\,dx-
\int_{\R^N}|\n \phi^{h_j}_{\mu}(x-a_j)|^2\,dx\bigg|\\
\nonumber&=
\int_{\mu^{-1}((\R^N\setminus\omega)-a_j)}
(1-\psi^2(\mu y+a_j))|\n \phi_{h_j}(y)|^2\,dy
=\begin{cases}
O\big(\mu^{N-2+2\sigma_{h_j}}\big),&\text{if }\mu_1(h_j)<N-1,\\[5pt]
O\big(\mu^{N-2\e}\big),&\text{if }\mu_1(h_j)\geq N-1,
\end{cases}
\end{align}
and
\begin{align}
\label{eq:68}
&\int_{\R^N}|\phi^{h_j}_{\mu}(x-a_j)|^2 |\n\psi(x)|^2\,dx\leq
{\rm const}\,\mu^2 \int_{\mu^{-1}((\Omega\setminus\omega)-a_j)}
|\phi_{h_j}(y)|^2\,dy\\
&\nonumber\quad
\leq {\rm const\,}
\mu^2
\int_{\mu^{-1}r}^{\mu^{-1}R}
s^{2(-\sigma_{h_j}-N+2)+N-1}\,ds=O(\mu^{2\sigma_{h_j}+N-2}),
\end{align}
where $r=\mathop{\rm dist}(a_j, \R^N\setminus\omega)$ and $R>0$ is such that 
$\Omega\subset B(a_j,R)$.
Similarly,
\begin{equation}\label{eq:69}
\int_{\R^N}\psi(x)\phi^{h_j}_{\mu}(x-a_j)\n\psi(x)\cdot \n
\phi^{h_j}_{\mu}(x-a_j)\, dx=O(\mu^{2\sigma_{h_j}+N-2})+
O\big(\mu^{N-2\e}\big).
\end{equation}
Estimate \eqref{eq:46} follows from (\ref{eq:67}--\ref{eq:69}).
The proof of (\ref{eq:48}--\ref{eq:51}) is analogous and is based on (\ref{eq:3}).
From
\begin{align*}
  S_{\Omega}(h_1,\dots,h_k) \notag\leq\  &
  \frac{{\displaystyle{\int_{\R^N}}}|\n
    (\psi(x)\phi^{h_j}_{\mu}(x-a_j))|^2\,dx-{\displaystyle{\int_{\R^N}}}
    \dfrac{h_j\big(\frac{x-a_j}
      {|x-a_j|}\big)}{|x-a_j|^2}|\psi(x)\phi^{h_j}_{\mu}(x-a_j)|^2
    \,dx} {\bigg({\displaystyle{\int_{\R^N}}}
    |\psi(x)\phi^{h_j}_{\mu}(x-a_j)|^{2^*}\,dx\bigg)^{\!\!2/2^*}
  }\\[5pt]
  & -
  \sum_{i\not=j}\frac{{\displaystyle{\int_{\R^N}}}\dfrac{h_i\big(\frac{x-a_i}
      {|x-a_i|}\big)}{|x-a_i|^2}|\psi(x)\phi^{h_j}_{\mu}(x-a_j)|^2\,dx
  } {\bigg({\displaystyle{\int_{\R^N}}}
    |\psi(x)\phi^{h_j}_{\mu}(x-a_j)|^{2^*}\,dx\bigg)^{\!\!2/2^*} },
\end{align*}
Lemmas \ref{l:interest1} and \ref{l:interest2}, 
and (\ref{eq:46}--\ref{eq:51}), it follows that
\begin{align*}
&S_{\Omega}(h_1,\dots,h_k)\leq 
S(h_j)+O(\mu^{2\sigma_{h_j}+N-2})+
O\big(\mu^{N-2\e}\big)\\
&-
\begin{cases}
  \mu^2\|\phi_{h_j}\|_{L^{2^*}(\R^N)}^{-2}\big(\int_{\R^N}\phi_{h_j}^2(x)\big)\Bigg(
  {\displaystyle{\sum_{i\neq j}}} \frac{h_i\big(\frac{a_j-a_i}
    {|a_j-a_i|}\big)}{|a_j-a_i|^2}+o(1)\Bigg),
  &\text{if }\mu_1(h_j)>-\big(\frac{N-2}2\big)^2+1,\\[18pt]
  \mu^2\|\phi_{h_j}\|_{L^{2^*}(\R^N)}^{-2}
  \big(\int_{|x|<\frac1\mu}\phi_{h_j}^2(x)\big)\Bigg(
  {\displaystyle{\sum_{i\neq j}}} \frac{h_i\big(\frac{a_j-a_i}
    {|a_j-a_i|}\big)}{|a_j-a_i|^2}+o(1)\Bigg), &\text{if
  }\mu_1(h_j)=-\big(\frac{N-2}2\big)^2+1,
\end{cases}
\end{align*}
as $\mu\to 0^+$. Since $0<\e<\frac{N-2}2$, there holds
$O\big(\mu^{N-2\e}\big)=o(\mu^2)$, thus implying the validity
of~(\ref{eq:47}).~\end{pf}

\begin{Corollary}\label{c:probl-bound-doma-1}
Let $j\in\{1,2,\dots,k\}$ such that $\mu_1(h_j)\geq-\big(\frac{N-2}2\big)^{2}+1$.
If  
\begin{align*}
\sum_{i\neq j}
\frac{h_i\big(\frac{a_j-a_i}
 {|a_j-a_i|}\big)}{|a_j-a_i|^2}>0,
\end{align*}
then
$$
S_{\Omega}(h_1,\dots,h_k)<S(h_j).
$$
\end{Corollary}
\begin{pf}
  It follows directly from Lemma \ref{lem:sub} after noticing that if
  $\mu_1(h_j)>-\big(\frac{N-2}2\big)^{2}+1$ then $2\sigma_{h_j}+N-2>2$
  and hence $O(\mu^{2\sigma_{h_j}+N-2})=o(\mu^{2})$ as $\mu\to 0^+$,
  while if $\mu_1(h_j)=-\big(\frac{N-2}2\big)^{2}+1$ then
  $2\sigma_{h_j}+N-2=2$ and hence
  $O(\mu^{2\sigma_{h_j}+N-2})=o\big(\mu^{2}\int_{|x|<\frac1\mu}
  \phi_{h_j}^2(x)\big)$. Taking $\mu$ sufficiently small, we obtain
  $S_{\Omega}(h_1,\dots,h_k)<S(h_j)$.
\end{pf}

\begin{pfn}{Theorem \ref{t:exlim}}
It follows from  Theorem \ref{th:psbound} and Corollary \ref{c:probl-bound-doma-1}
arguing as in the proof of Theorem \ref{t:ach}.
\end{pfn}

\end{document}